\RequirePackage[reqno]{amsmath}
\documentclass[preprint]{amsart}

\usepackage{amssymb}
\usepackage{empheq}
\usepackage{amsthm}
\usepackage{enumerate}
\usepackage{comment}
\usepackage[top=1in, left=1in, right=0.9in, bottom=1in]{geometry}
\usepackage{mathrsfs}

\newtheorem{theorem}{Theorem}[section]
\newtheorem{lemma}[theorem]{Lemma}
\newtheorem{definition}[theorem]{Definition}
\newtheorem{proposition}[theorem]{Proposition}
\newtheorem{corollary}[theorem]{Corollary}

\def\D{{\mathbb D}}

\def\R{{\mathbb R}}

\def\S{{\mathbb S}}

\newcommand{\cF}{\mathcal{F}}

\newcommand{\ds}{\displaystyle{}}

\newcommand{\overbar}[1]{\mkern 1.5mu\overline{\mkern-1.5mu#1\mkern-1.5mu}\mkern 1.5mu}

\newcommand{\cp}{\operatorname{cap}}
\newcommand{\supp}{\operatorname{supp}}

\newcommand{\Beta}{{\rm B}}

\begin{document}

\title{Minimum Riesz Energy Problem on the Hyperdisk}

\author[Mykhailo Bilogliadov]{Mykhailo Bilogliadov}
\email{mykhail@okstate.edu}
\address{Department of Mathematics, Oklahoma State University, Stillwater, OK 74078, U.S.A.}


\begin{abstract}
We consider the minimum Riesz $s$-energy problem on the unit disk $\D:=\{(x_1,\ldots,x_d)\in\R^d: x_1=0, x_2^2+x_3^2+\ldots+x_d^2\leq 1\}$ in the Euclidean space $\R^d$, $d\geq 3$, immersed into a smooth rotationally invariant external field $Q$. The charges are assumed to interact via the Riesz potential $1/r^s$, with $d-3 < s < d-1$, where $r$ denotes the Euclidean distance. We solve the problem by finding an explicit expression for the extremal measure. We then consider applications to a monomial external field and an external field generated by a positive point charge, located at some distance above the disk on the polar axis. We obtain an equation describing the critical height for the location of the point charge, which guarantees that the support of the extremal measure occupies the whole disk $\D$. We also show that under some mild restrictions on a general external field the support of the extremal measure will have a ring structure. Furthermore, we demonstrate how to reduce the problem of recovery of the extremal measure in this case to a Fredholm integral equation of the second kind.

\smallskip
\noindent \textbf{MSC 2010.} 31B05, 31B10, 31B15, 45B05, 45D05

\smallskip
\noindent \textbf{Key words.} Minimal energy problem, Riesz potential, weighted energy, equilibrium measure, extremal measure
\end{abstract}

\maketitle

\medskip
\medskip
\noindent

\section{Introduction}

Let $\S^{d-1}:=\{x\in\R^d: |x|=1\}$ be the unit sphere in $\R^d$, and $\D_R:=\{(x_1,\ldots,x_d)\in\R^d: x_1=0, x_2^2+x_3^2+\ldots+x_d^2\leq R^2\}$ be the disk of radius $R$ in $\R^d$, with $d\geq 3$, and where $|\cdot|$ is the Euclidean distance. The ring $\mathscr{R}(a,b)$  in $\R^d$ is defined as $\mathscr{R}(a,b):=\{(0,r\overbar{x})\in\R^d: a\leq r\leq b, \overbar{x}\in\S^{d-2} \}$, and the unit disk in $\R^d$ will be denoted by $\D$. Given a compact set $E\subset\D$, consider the class $\mathcal{M}(E)$ of unit positive Borel measures supported on $E$. For $0<s<d$, the Riesz $s$-potential and Riesz $s$-energy of a measure $\mu\in\mathcal{M}(E)$ are defined respectively as 
$$U_s^{\mu}(x):=\int \frac {1} {|x-y|^s} \, d\mu(y), \quad I_s(\mu):=\iint \frac {1} {|x-y|^s} \, d\mu(x)d\mu(y).$$
Let 
$$W_s(E):=\inf\{I_s(\mu): \mu\in\mathcal{M}(E)\}.$$
Define the Riesz $s$-capacity of $E$ as $\cp_s(E):=1/W_s(E)$. We say that a property holds quasi-everywhere (q.e.), if the exceptional set has a Riesz $s$-capacity zero. When $\cp_s(E)>0$, there is a unique $\mu_{E}$ such that $I_s(\mu_E)=W_s(E)$. Such $\mu_E$ is called the Riesz $s$-equilibrium measure for $E$.

An external field is defined as a non-negative lower-semicontinuous function $Q: E \rightarrow [0,\infty]$, such that $Q(x)<\infty$ on a set of positive surface area measure. The weighted energy associated with $Q(x)$ is then defined by $$I_Q(\mu):=I_s(\mu)+2\int Q(x)d\mu(x).$$

\noindent The minimum energy problem on $\D$ in the presence of the external field $Q(x)$ refers to the minimal quantity 
$$V_Q:=\inf\{I_Q(\mu): \mu\in\mathcal{M}(E)\}.$$ 
A measure $\mu_Q \in\mathcal{M}(E)$ such that $I_Q(\mu_Q)=V_Q$ is called the $s$-extremal (or positive Riesz $s$-equilibrium) measure associated with $Q(x)$.

The potential $U_s^{\mu_Q}$ of the measure $\mu_Q$ satisfies the Gauss variational inequalities
\begin{empheq}{align}\label{var1}
U_s^{\mu_Q}(x)+Q(x)&\geq F_Q \quad \text{q.e. on}\ \ E,\\\label{var2}
U_s^{\mu_Q}(x)+Q(x)& \leq F_Q \quad \text{for all}\ \ x\in S_Q,
\end{empheq}
where $F_Q:=V_Q-\ds{ \int Q(x)\,d\mu_Q(x)}$, and $S_Q:=\supp{\mu_Q}$, which is sometimes refferred as the extremal support (see Theorem 10.3 in \cite{mizuta} or Proposition 3 in \cite{bds2}, and also book \cite{bhs}). We remark that for continuous external fields, the inequality in $(\ref{var1})$ holds everywhere, which implies that equality holds in $(\ref{var2})$.

The minimum energy problems with external fields on the sphere and other manifolds have been a subject of intensive study by the group of 
Brauchart, Dragnev and Saff, see \cite{bds1}--\cite{bds3}, \cite{ddss}; see also \cite{bilo}. For a comprehensive treatment of the 
subject one should refer to a forthcoming book \cite{bhs}.

The main aim of this paper is to provide a solution to the weighted energy problem on the unit disk 
$\D\subset\R^d$, $d\geq 3$, immersed in a general smooth external field, possessing rotational symmetry with respect to the polar 
axis. It is assumed that the charges interact according to the Riesz $s$-potential $1/r^s$, with $d-3<s<d-1$, where $r$ is the Euclidean distance. We obtain an explicit expression for the density of the $s$-extremal measure, assuming that the support of the $s$-extremal measure occupies the whole disk $\D$. Our results generalize the solution of the problem obtained by Copson \cite{cop}, who considered this problem for the classical Coulomb potential in $\R^3$.

The first application of our results is concerned with a situation when the disk $\D$ is immersed into a monomial external field. We find the extremal measure corresponding to such a field, while also explicitly finding the extremal support. Our second application is concerned with finding the extremal measure when an external field is generated by a positive point charge, placed on the polar axis at a certain distance above the disk $\D$. A similar problem for the sphere $\S^2$ for the Coulomb potential in $\R^3$ was raised by Gonchar \cite{lmfns}, and solved in \cite{bds1} for general Riesz potentials on the sphere $\S^{d-1}\subset\R^d$, $d\geq 3$. The extensions and further results on Gonchar's problem are contained in works \cite{bds3}, \cite{bdsw}. The problem of finding a signed measure representing the charge distribution on the disk $\D$ in $\R^3$ under the influence of a positive point charge placed on the polar axis above the disk $\D$, for the case of Coulomb potential, was first considered in the classical work of Thomson \cite{thomson}, and later solved by a different method by Gallop \cite{gallop}. Below we solve this problem for the case of higher dimensions and general Riesz $s$-potentials by finding the extremal measure representing the positive charge distribution on the disk $\D$.  For the case of higher dimensions and general Riesz $s$-potentials we give an explicit estimate on the height and magnitude of the point charge, which guarantees for the extremal support to occupy the whole disk $\D$. Moreover, in the case of classical Coulomb potential in $\R^3$ and a positive unit point charge, we are able to precisely determine the height of the point charge so that the extremal support occupies the entire disk $\D$.

Moreover, we investigate what happens to the support of the extremal measure $\mu_Q$ if one moves a point charge closer to the disk, beyond the aforementioned height threshold. It is demonstrated that under some mild restrictions on a general external field, the support $\supp\mu_Q$ will be a ring, contained in the disk $\D$. We derive an integral equation allowing to recover the extremal measure $\mu_Q$, when $\supp\mu_Q$ is a ring.

We start by taking advantage of the rotationally symmetry of $\D$ using the cylindrical coordinates $z,r, \theta_1, \theta_2$, $\ldots, \theta_{d-3}, \varphi$, defined as
\begin{align*}
x_1 & = z,\\
x_2 & = r\cos\theta_1,\\
x_3 & = r\sin\theta_1\cos\theta_2,\\
x_3 & = r\sin\theta_1\sin\theta_2\cos\theta_3,\\
& \vdots \\
x_{d-2} & = r\sin\theta_1\sin\theta_2\ldots\sin\theta_{d-4} \cos\theta_{d-3},\\
x_{d-1} & = r\sin\theta_1\sin\theta_2\ldots\sin\theta_{d-3} \cos\varphi,\\
x_{d} & = r\sin\theta_1\sin\theta_2\ldots\sin\theta_{d-3} \sin\varphi,
\end{align*}
where $r\geq0$, $0\leq\theta_j\leq \pi$, $j=1,2,\ldots,d-3$, and $0\leq\varphi \leq 2\pi$. The surface area element on a surface of constant height $z$, written in cylindrical coordinates, is given by
\[
d S = r^{d-2}\, \sin^{d-3} \theta_1 \, \sin^{d-4} \theta_2 \ldots \sin\theta_{d-3} \, dr\, d\theta_1\, d\theta_2 \ldots d\theta_{d-3}\, d\varphi = r^{d-2}\, dr\, d\sigma_{d-1},
\]
where $\sigma_d$ is the surface area element of the unit sphere $\S^{d-1}$. The total surface area of the sphere $\S^{d-1}$ is given by 
\[
\omega_d=\frac{2\pi^{d/2}}{\Gamma(d/2)}.
\]

\noindent In what follows, we will need to use certain special functions, for which we fix the notation here. The incomplete Beta function $\Beta(z; a,b)$ is defined as
\begin{equation}\label{betafdef}
\Beta(z;a,b) := \int_0^z t^{a-1} (1-t)^{b-1}\, dt,
\end{equation}
and the Beta function $\Beta(a,b):=\Beta(1; a,b)$. The Gauss hypergeometric function $_2 F_1(a,b;c,z)$ is defined via series
\begin{equation}\label{gausshyperdef}
_2 F_1(a,b;c,z) := \sum_{n=0}^\infty \frac {(a)_n\, (b)_n} {(c)_n} \, \frac {z^n} {n!}, \quad |z|<1,
\end{equation}
where $(a)_0: = 1$ and $(a)_n: = a (a+1) \ldots (a+n-1)$ for $n\geq 1$ is the Pochhammer symbol.

We commence by recording the sufficient conditions on an external field $Q$ that guarantee that the support of the extremal measure $\mu_Q$ is a ring or a disk.

\begin{theorem}\label{theo0}
Let $s=(d-3)+2\lambda$, with $0<\lambda<1$. Assume that the external field $Q: \D\to [0,\infty]$ is invariant with respect to the rotations about the polar axis, that is $Q(x)=Q(r)$, where $x=(0,r\overbar{x})\in\D$, $\overbar{x}\in\S^{d-2}$, $0\leq r \leq 1$. Further suppose that $Q$ is a convex function, that is $Q(r)$ is convex on $[0,1]$. Then the support of the extremal measure $\mu_Q$ is a ring $\mathscr{R}(a,b)$, contained in the disk $\D$. In other words, there exist real numbers $a$ and $b$ such that $0\leq a < b\leq 1$, so that $\supp\mu_Q=\mathscr{R}(a,b)$.

Furthermore, if $Q(r)$ is, in addition, an increasing function, then $a=0$, which implies that the support of the extremal measure $\mu_Q$ is a disk of radius $b\leq1$, centered at the origin. 

On the other hand, if $Q(r)$ is a decreasing function, then $b=1$, that is the support of the extremal measure $\mu_Q$ will be a ring with outer radius $1$. 
\end{theorem}

The support $S_Q$ is a main ingredient in determining the extremal measure $\mu_Q$ itself. Indeed, if 
$S_Q$ is known, the equilibrium measure $\mu_Q$ can be recovered by solving the singular integral equation
\begin{equation}\label{ieec}
\int \frac{1}{|x-y|^s}\,d\mu(y)+Q(x)=F_Q, \quad x\in S_Q,
\end{equation}
where $F_Q$ is a constant (see (\ref{var2})).

We solve this equation and obtain the following theorem, which explicitly gives the density of the extremal measure when the support $S_Q$ is the disk $\D_R$. Our results extend the original work of Copson \cite{cop}, which dealt with classical Coulomb potential in $\R^3$.
\begin{theorem}\label{theo1}
Suppose that the support of the extremal measure $\mu_Q$ is the disk $\D_R$, and the external field $Q$ is invariant with respect to rotations about the polar axis, that is $Q(x)=Q(r)$, where $x=(0,r\overbar{x})\in\D_R$, $\overbar{x}\in\S^{d-2}$, $0\leq r \leq R$. Also assume that $Q\in C^2(\D_R)$.  Let $s=(d-3)+2\lambda$, with $0<\lambda<1$, and let
\begin{equation}\label{func-F}
F(t) = \frac {\sin(\lambda\pi)\, \Gamma((d-3)/2+\lambda)}{\pi^{(d+1)/2}\, \Gamma(\lambda)}\, \frac{1}{t}\, \frac{d}{dt} \int_t^R \frac {g(r)\, r\, dr} {(r^2-t^2)^{1-\lambda}}, \quad 0\leq t \leq R,
\end{equation}
with
\begin{equation}\label{func-g}
g(r) = \frac {1}{r^{d+2\lambda-4}}\, \frac{d}{dr} \int_0^r \frac {Q(u)\, u^{d-2}\, du} {(r^2-u^2)^{1-\lambda}}, \quad 0\leq r \leq R.
\end{equation}
Then for the extremal measure $\mu_Q$ we have
\begin{equation}\label{meas-abs-cont}
d\mu_Q(x) = f(r)\, r^{d-2}\, dr\, d\sigma_{d-1}(\overbar{x}), \quad x=(0,r\overbar{x})\in\D_R, \quad \overbar{x}\in\S^{d-2}, \quad 0 \leq r \leq R,
\end{equation}
where the density $f$ is explicitly given by
\begin{equation}\label{extremal-density}
f(r) = C_Q\, (R^2-r^2)^{\lambda-1} + F(r), \quad 0\leq r \leq R,
\end{equation}
with the constant $C_Q$ uniquely defined by
\begin{equation}\label{constant-C}
C_Q = \frac {2\, \Gamma((d-1)/2+\lambda)} {\Gamma(\lambda)\, \Gamma((d-1)/2)}\, \frac{1}{R^{d+2\lambda-3}} \, \bigg\{ \frac {\Gamma((d-1)/2)} {2\, \pi^{(d-1)/2}} - \int_0^R F(t)\, t^{d-2}\, dt \bigg\}.
\end{equation}

\end{theorem}

\section{Applications}

In what follows we will need to know the equilibrium measure and capacity of a disk $\D_R$ of radius $R>0$. The following theorem extends the corresponding result of Copson \cite{cop}, which dealt with the Coulomb potential in $\R^3$.

\begin{theorem}\label{equil-meas-cap-disk}
Let $s=(d-3)+2\lambda$, with $0<\lambda<1$. The equilibrium measure $\mu_{\D_R}$ of the disk $\D_R$ of radius $R$ is given by
\begin{equation}\label{equil-meas-disk}
d\mu_{\D_R}(x) = f(r)\, r^{d-2}\, dr\, d\sigma_{d-1}(\overbar{x}), \quad x=(0,r\overbar{x})\in\D_R, \quad \overbar{x}\in\S^{d-2}, \quad 0 \leq r \leq R,
\end{equation}
where the density $f$ is 
\begin{equation}\label{equil-density-disk}
f(r) =\frac {\Gamma((d+2\lambda-1)/2)} {\pi^{(d-1)/2}\, \Gamma(\lambda)\, }\, \frac{1}{R^{d+2\lambda-3}}\, (R^2-r^2)^{\lambda-1}, \quad 0\leq r \leq R.
\end{equation}
The capacity of the disk $\D_R$ is given by
\begin{equation}\label{cap-disk}
\cp_s (\D_R) = \frac  {\sin(\lambda\pi) \Gamma(\lambda) \Gamma((d-1)/2)} {\pi\Gamma((d+2\lambda-1)/2)}\, R^{d+2\lambda-3}.
\end{equation}
\end{theorem}

Assume that the disk $\D$ is immersed into a general rotationally invariant external field $Q$, satisfying the conditions of the second statement of Theorem $\ref{theo0}$. It then follows that the support of the extremal measure $\mu_Q$ will be a disk $\D_R$ of some radius $R\leq1$. The (presently) unknown radius $R$ can be found by minimizing the Mhaskar-Saff functional, which is defined as follows \cite{bds1}.
\begin{definition}\label{ms-func-def}
The  $\cF$-functional of a compact subset $E\subset \D$ of positive Riesz $s$-capacity is defined as
\begin{equation}\label{ms-func}
\cF_s(E) := W_s(E) + \int Q(x)\,d\mu_E(x), 
\end{equation}
where $W_s(E)$ is the Riesz $s$-energy of the compact $E$ and $\mu_E$ is the equilibrium measure (with no external field) on $E$.
\end{definition}
The main objective of introducing the $\cF_s$-functional is its following extremal property, which originally was proved in \cite{bds1} for the general Riesz potentials on the sphere $\S^{d-1}$, $d\geq 3$.
\begin{proposition}\label{ms-func-min-prop}
Let $Q$ be an external field on $\D$. Then $\cF_s$-functional is minimized for $S_Q=\supp(\mu_Q)$. 
\end{proposition}
Utilizing Proposition \ref{ms-func-min-prop}, we can now explicitly determine the support of the extremal measure provided the external field satisfies some mild restrictions.  
\begin{theorem}\label{critical-radius-theo}
Let $s=(d-3)+2\lambda$, with $0<\lambda<1$. Assume that the external field $Q: \D\to [0,\infty]$ is invariant with respect to the rotations about the polar axis, that is $Q(x)=Q(r)$, where $x=(0,r\overbar{x})\in\D$, $\overbar{x}\in\S^{d-2}$, $0\leq r \leq 1$. Assume that $Q\in C^2(\D)$. Further suppose that $Q$ is a convex increasing function, that is $Q(r)$ is convex increasing on $[0,1]$. Then the support of the extremal measure $\mu_Q$ will be a disk of radius $R\leq1$, centered at the origin. The radius $R$ of this disk is either the unique solution of the equation 
\begin{equation}\label{radius-equation-equation}
\frac{2\sin(\lambda\pi)}{\pi(d+2\lambda-3)} \int_0^R Q'(r)\, (R^2-t^2)^{\lambda-1}\, t^{d-1}\, dt =1,
\end{equation}
on the interval $(0,1]$ if it exists, or $R=1$ when such a solution fails to exist.
\end{theorem}

As the first applications of our results, we consider the situation when the disk $\D$ is immersed into an external field given by a monomial, namely
\begin{equation}\label{ext-field-monom-gen}
Q(x) = q r^\alpha, \quad q > 0, \quad \alpha \geq 1, \quad  x=(0,r\overbar{x})\in\D, \quad \overbar{x}\in\S^{d-2}, \quad 0 \leq r \leq 1.
\end{equation}

\noindent It is clear that external field $Q$ in $(\ref{ext-field-monom-gen})$ is invariant with respect to the rotations about the polar axis. Also, $Q(r)$ is a non-negative increasing convex function on $[0,1]$. From Theorem \ref{theo0} it then follows that the support of the corresponding extremal measure $\mu_Q$ will be a disk $\D_R$, with some $R\leq1$. First invoking Theorem \ref{critical-radius-theo} we compute the extremal support, and then with that knowledge at hand we use Theorem \ref{theo1} to find a closed-form expression for the extremal measure. 
\begin{theorem}\label{mon-ext-field-extrem-meas}
Let $s=(d-3)+2\lambda$, with $0<\lambda<1$. The extremal measure $\mu_Q$, corresponding to the monomial external field $(\ref{ext-field-monom-gen})$, is supported on the disk $\D_{R_\ast}$, where $R_\ast$ is defined as
\begin{equation}\label{mon-ext-field-extrem-meas-disk-rad}
R_\ast = \left(\frac{(d+2\lambda-3)\pi\Gamma((d+\alpha+2\lambda-1)/2)}{q\alpha\sin(\lambda\pi)\Gamma(\lambda)\Gamma((d+\alpha-1)/2)}\right)^{1/(d+\alpha+2\lambda-3)}.
\end{equation}
For the extremal measure $\mu_Q$ we have
\begin{equation}\label{mon-ext-field-extrem-meas-formula}
d\mu_Q(x) = f(r)\, r^{d-2}\, dr\, d\sigma_{d-1}(\overbar{x}), \quad x=(0,r\overbar{x})\in\D_{R_\ast}, \quad \overbar{x}\in\S^{d-2}, \quad 0 \leq r \leq R_\ast,
\end{equation}
with the density $f(r)$ is given by
\begin{equation}\label{mon-ext-field-extrem-dens}
f(r) = C_Q\, (R^2-r^2)^{\lambda-1} + F(r), \quad 0 \leq r \leq R_\ast,
\end{equation}
where
\begin{empheq}{align}\label{mon-ext-field-func-F}
 F(r) & =  \frac{q\,\sin(\lambda\pi)\,\Gamma((d+\alpha-1)/2)\, \Gamma((d+2\lambda-3)/2)}{\pi^{(d+1)/2}\, \Gamma((d+\alpha+2\lambda-3)/2)}\, R_\ast^\alpha\, (R_\ast^2-r^2)^{\lambda-1}\times  \\\nonumber
 \bigg\{ & - {}_2 F_1\left(-\frac{\alpha}{2}, 1;\lambda+1; 1-\left(\frac{r}{R_\ast}\right)^2\right) \\\nonumber
 & + \frac{\alpha}{2\lambda(\lambda+1)}\bigg(1-\left(\frac{r}{R_\ast}\right)^2 \bigg)\, {}_2 F_1\left(1-\frac{\alpha}{2}, 2;\lambda+2; 1-\left(\frac{r}{R_\ast}\right)^2\right) \bigg\}, \quad   0\leq r\leq R_\ast,
\end{empheq}
and the constant $C_Q$ is defined as
\begin{equation}\label{mon-ext-field-func-C-Q}
C_Q =   \frac {\Gamma((d+2\lambda-1)/2)} {\pi^{(d-1)/2}\, \Gamma(\lambda)} \, \bigg\{ \frac{1}{R_\ast^{d+2\lambda-3}} +  \frac{q\sin(\lambda\pi)\,\Gamma((d+\alpha-1)/2)\, \Gamma(\lambda)}{\pi\, \Gamma((d+\alpha+2\lambda-1)/2)}\, R_\ast^\alpha \bigg\}.
\end{equation}

\end{theorem}


Another application is concerned with finding the extremal measure $\mu_Q$ in the case of the Riesz $s$-potential generated by a positive point charge. We assume that the external field $Q$ is produced by a positive point charge of magnitude $q$, placed on the positive polar semi-axis at some distance $h>0$ above the disk $\D$. This external field is given by
\begin{equation}\label{pt-charge}
Q(x)=\frac{q}{(r^2+h^2)^{s/2}}, \quad h > 0, \quad q>0, \quad s>0, \quad x=(0,r\overbar{x})\in\D, \quad \overbar{x}\in\S^{d-2}, \quad 0 \leq r \leq 1.
\end{equation}
This problem is similar to the celebrated Gonchar's problem, which was solved in \cite{bds1} for the case of classical Newtonian potential in $\R^d$. The Gonchar's problem is concerned with a situation when a positive unit point charge is approaching the insulated unit sphere, carrying a total charge $1$, eventually causing a spherical cap free of charge to appear. Gonchar raised a question about finding the smallest distance from the point charge to the sphere such that the whole of the sphere still being positively charged. In a slightly more general setting, the solution of the Gonchar's problem means that if a point charge of a non-negative magnitude $q$, located on the positive polar semi-axis, is too far from the surface of the sphere $\S^2$, or the magnitude $q$ of the point charge is too small, the electrostatic field, created by this point charge, is too weak to force the equilibrium charge distribution from occupying the whole surface of the sphere $\S^2$. It is known that in the case of the positive unit point charge, the critical height is precisely $1+\rho$, where $\rho$ is the golden ratio $(1+\sqrt{5})/2  \approx 1.6180339887$.

The question about the charge distribution on the surface of the disk $\D$ in $\R^3$, influenced by a positive point charge placed above the disk on the polar axis, was first considered by Thompson \cite{thomson}. By a different method, the problem was treated by Gallop \cite{gallop}.


\begin{theorem}\label{theo3} Let $s=(d-3)+2\lambda$, with $0<\lambda<1$. Assume that the external field $Q$ is given by $(\ref{pt-charge})$, with $h>\max\{h_-,h_+\}$, where
\begin{equation}\label{h-minus}
h_-:=\left(\frac{q\, ((1-\lambda)(d-2\lambda+1)+1)^2 \sin(\lambda\pi)}{8\pi(d+2\lambda-1)(1-\lambda)} \, \Beta\left(\lambda, \frac{d-1}{2}\right) \right)^{1/(d+2\lambda-3)},
\end{equation}
and $h_+$ is the largest positive root of the function
\begin{empheq}{align}\label{func-p}
p(h) & = \frac {2\, \Gamma((d-1)/2+\lambda)} {\Gamma(\lambda)\, \Gamma((d-1)/2)}\, \bigg\{ \frac {\Gamma((d-1)/2)} {2\, \pi^{(d-1)/2}} + q\, c_{d,\lambda} \bigg\}  \\\nonumber
& - \frac {q\,\sin(\lambda\pi)\,\Gamma((d-1)/2)} {\pi^{(d+1)/2}} \,  \bigg\{ \frac{d-2\lambda-1}{2}\,  \frac{1} {h^{d-1}}\, \Beta\left( \frac{1}{1+h^2}; \lambda, \frac{d-2\lambda-1}{2}\right)  + \frac{1} {h^{2\lambda} (1+h^2)^{(d-3)/2}} \bigg\},
\end{empheq}
with
\begin{empheq}{align}\label{pt-charge-small-c}
c_{d,\lambda}:= \frac{\sin(\lambda\pi)\,\Gamma((d-1)/2)}{\pi^{(d+1)/2}}\, & h^{2(1-\lambda)} \,  \bigg\{ \frac{d-2\lambda-1}{2}\, \int_0^1 \frac{t^{d-2}}{(h^2+t^2)^{(d-2\lambda+1)/2}}\, \Beta\left( \frac{1-t^2}{1+h^2}; \lambda, \frac{d-2\lambda-1}{2}\right) dt \\ \nonumber
& + \frac {\Gamma(\lambda)\, \Gamma((d-1)/2)}  {2\, \Gamma((d-1)/2+\lambda)}  \, \frac{1}{(1+h^2)^{(d-1)/2}}\, {}_2 F_1\left(1, \lambda; \frac{d+2\lambda-1}{2}; \frac{1}{1+h^2}\right) \bigg\}.
\end{empheq}
The extremal measure $\mu_Q$, corresponding to the external field of a point charge $(\ref{pt-charge})$, is given by
\begin{equation}\label{pt-charge-ext-field-ext-mes}
d\mu_Q(x) = f(r)\, r^{d-2}\, dr\, d\sigma_{d-1}(\overbar{x}), \quad x=(0,r\overbar{x})\in\D, \quad \overbar{x}\in\S^{d-2}, \quad 0 \leq r \leq 1,
\end{equation}
where
\begin{equation}\label{pt-charge-ext-field-dens}
f(r) = C_Q\, (1-r^2)^{\lambda-1} + F(r), \quad 0\leq r \leq 1.
\end{equation}
Here 
\begin{empheq}{align}\label{pt-charge-func-F-theo}
F(r) = - \frac{q\,\sin(\lambda\pi)\,\Gamma((d-1)/2)}{\pi^{(d+1)/2}}\, h^{2(1-\lambda)} \, & \bigg\{ \frac{d-2\lambda-1}{2}\,  \frac{1}{(h^2+r^2)^{(d-2\lambda+1)/2}}\, \Beta\left( \frac{1-r^2}{1+h^2}; \lambda, \frac{d-2\lambda-1}{2}\right) \\ \nonumber
& + \frac{(1-r^2)^{\lambda-1}} {(1+h^2)^{(d-3)/2}\, (h^2+r^2)} \bigg\}, \quad 0\leq r \leq 1,
\end{empheq}
and where the positive constant $C_Q$ is given by
\begin{equation}\label{pt-charge-C_Q}
C_Q = \frac {2\, \Gamma((d-1)/2+\lambda)} {\Gamma(\lambda)\, \Gamma((d-1)/2)}\, \bigg\{ \frac {\Gamma((d-1)/2)} {2\, \pi^{(d-1)/2}} + q\, c_{d,\lambda} \bigg\},
\end{equation}
\end{theorem}

\noindent In the important special case of Newtonian potential in the even dimensions starting with $8$, from the above Theorem it follows that the extremal measure can be written in a simplified form.
\begin{corollary}\label{pnt-charge-newtonian-case}
 Let $d=2m+4$ and $s=d-2=2(m+1)$, where $m\geq 2$. Assume that the external field $Q$ is given by $(\ref{pt-charge})$, with $h>\max\{h_-,h_+\}$, with
\begin{equation}\label{h-minus-newt-pot}
h_-:=\left(\frac{q\, (m+2)\,\Gamma(m+3/2)}{8\sqrt{\pi}(m+1)!} \right)^{1/2(m+1)},
\end{equation}
and $h_+$ is the largest positive root of the function
\begin{empheq}{align}\label{func-p-newt-pot}
p(h) & = \frac {2\, (m+1)!} {\sqrt{\pi}\, \Gamma(m+3/2)}\, \bigg\{ \frac { \Gamma(m+3/2)} {2\, \pi^{m+3/2}} + q c  \bigg\}  \\\nonumber
& - \frac {q\, \Gamma(m+3/2)} {\pi^{m+5/2}} \,  \bigg\{ 2(m+1) \, h^{-(2m+3)} \, \sum_{n=0}^m \frac{(-m)_n}{(2n+1)n!}\,(1+h^2)^{-(n+1/2)}  + \frac{1} {h (1+h^2)^{m+1/2}} \bigg\},
\end{empheq}
with
\begin{empheq}{align}\label{pt-charge-small-c-newt-pot}
c: &=  \frac{(\Gamma(m+3/2))^2}{\pi^{m+5/2}}\, h\, (1+h^2)^{-(m+1)} \\ \nonumber
&  \bigg\{  \Gamma\left(m+\frac{3}{2}\right)\, (1+h^2)^{-(m+5/2)}\, \sum_{n=0}^m \frac{(-m)_n}{(2n+1)n!}\,(1+h^2)^{-n} \, \sum_{l=0}^{m-2} \frac{(2-m)_l\,\Gamma(n+l+3/2)}{(n+m+l+2)!}\,(1+h^2)^{-l} \\ \nonumber
& + \frac{\sqrt{1+h^2}}{h+\sqrt{1+h^2}}\, \sum_{n=0}^m \frac{(-m)_n}{(m+n+1)!}\, \left(\frac{\sqrt{1+h^2}-h}{\sqrt{1+h^2}+h}\right)^n \bigg\}.
\end{empheq}
Then the extremal measure $\mu_Q$, corresponding to the external field of a point charge $(\ref{pt-charge})$, is given by
\begin{equation}\label{pt-charge-ext-field-ext-mes-newt-pot}
d\mu_Q(x) = f(r)\, r^{d-2}\, dr\, d\sigma_{d-1}(\overbar{x}), \quad x=(0,r\overbar{x})\in\D, \quad \overbar{x}\in\S^{d-2}, \quad 0 \leq r \leq 1,
\end{equation}
where
\begin{equation}\label{pt-charge-ext-field-dens-newt-pot}
f(r) = C_Q\, \frac{1}{\sqrt{1-r^2}} + F(r), \quad 0\leq r \leq 1.
\end{equation}
Here 
\begin{empheq}{align}\label{pt-charge-func-F-newt-pot}
F(r)= - qh \frac{\Gamma(m+3/2)}{\pi^{m+5/2}}\,\bigg\{ & \frac{ 2(m+1)}{(h^2+r^2)^{m+2}}\,   \sum_{n=0}^m \frac{(-m)_n}{(2n+1)n!}\, \left( \frac{1-r^2}{1+h^2}\right)^{n+1/2} \\ \nonumber
& + \frac{1}{\sqrt{1-r^2}}\, \frac{1}{(h^2+r^2)\,(1+h^2)^{m+1/2}} \bigg\}, \quad 0\leq r\leq 1.
\end{empheq}
and where the positive constant $C_Q$ is given by
\begin{equation}\label{pt-charge-C_Q-newt-pot}
C_Q =  \frac{2(m+1)!}{\sqrt{\pi}\Gamma(m+3/2)}\,\bigg\{ \frac{\Gamma(m+3/2)}{2\pi^{m+3/2}} + qc  \bigg\}.
\end{equation}
\end{corollary}

The case of the three-dimensional Euclidean space $\R^3$ and Coulomb potential, corresponding to $d=3$ and $\lambda=1/2$ in the context of Theorem \ref{theo3}, deserves special attention. Assuming that the disk $\D$ is immersed into the external field generated by a positive unit point charge, in this physically important case we are able to precisely determine the height of the point charge that guarantees the extremal support $\mu_Q$ to occupy the whole disk $\D$. This is an improvement of Theorem \ref{theo3}, where we can only provide an estimate of such a height.

\begin{corollary}\label{disk-three-dim-coulomb-pot}
Suppose the external field $Q$ is given by $(\ref{pt-charge})$, with $d=3$ and $s=1$, and where $h$ is chosen such that $h\geq h_+$, where $h_+$ is the unique positive root of the function 
\begin{equation*}
p(h)=\frac{1}{2\pi}\left(1+\frac{2h\tan^{-1} (1/h)}{\pi\sqrt{1+h^2}} \right) - \frac{1}{\pi^2 h} - \frac{1}{\pi^2h^2}\tan^{-1} (1/h).
\end{equation*}
Then, under these assumptions $S_Q=\D$, and the extremal measure $\mu_Q$ is given by
\begin{equation}\label{mes-disk-three-dim-coulomb-pot}
d\mu_Q(x) = f(r)\, r\, dr\, d\sigma_{2}(\overbar{x}), \quad x=(0,r\overbar{x})\in\D, \quad \overbar{x}\in\S^1, \quad 0 \leq r \leq 1,
\end{equation}
where the density $f(r)$ is
\begin{empheq}{align}\label{density-disk-three-dim-coulomb-pot}
 f(r) =  \frac{1}{2\pi}\left(1+\frac{2h\tan^{-1} (1/h)}{\pi\sqrt{1+h^2}} \right)  \frac{1}{\sqrt{1-r^2}} & - \frac{h}{\pi^2(h^2+r^2)}\frac{1}{\sqrt{1-r^2}}\\\nonumber
	& -\frac{h}{\pi^2}\frac{1}{(h^2+r^2)^{3/2}}\tan^{-1} \sqrt{\frac{1-r^2}{h^2+r^2}},\quad 0\leq r \leq 1.
\end{empheq}
If the height of the point charge is chosen such that $h < h_+$, then the support of the extremal measure $\mu_Q$ will no longer be the entire disk $\D$, as there will be an opening around the origin.
\end{corollary}

Note that from Corollary \ref{disk-three-dim-coulomb-pot} it follows that if the charge is moved closer to the disk past the critical height $h_+$, then support of the extremal measure will no longer be the entire disk. This means that when $h<h_+$ the point charge clears out an opening in the disk $\D$ at the origin, which will be free of charge and is likely to have a ring structure.

We now turn to the question of recovering the extremal measure $\mu_Q$ in the case when its support is a ring $\mathscr{R}(a,b)$, that is $\supp\mu_Q=\mathscr{R}(a,b)$.

Apparently the first attempt to solve integral equation $(\ref{ieec})$ on a ring was undertaken by Gubenko and Mossakovskii in 1960 for the case of classical Coulomb potential in $\R^3$. In paper \cite{gubmos} they considered a problem of calculating the pressure that a rigid die, having a form of a circular concentric ring, exerts on an elastic half-space. They were able to obtain an approximate solution to $(\ref{ieec})$, with a prescribed degree of accuracy, acceptable for their needs. The next major step was undertaken by Cooke in 1963, while seeking to obtain a closed form solution to equation $(\ref{ieec})$ on a ring, again for the case of classical Coulomb potential in $\R^3$. In paper \cite{cooke} he found a way to reduce equation $(\ref{ieec})$ to a Fredholm integral equation of the second kind. The approach of Cooke was based on exploiting certain representations of the kernel of integral equation $(\ref{ieec})$, based on intricate identities for Bessel functions. Williams \cite{williams} noticed that Cooke's result is, in fact, independent of any special functions. A further step in investigating equation $(\ref{ieec})$ on a ring was taken by Clements and Love in 1970s. In papers \cite{loveclem} and \cite{loveclem1}, instead of reducing $(\ref{ieec})$ to one Fredholm integral equation of the second kind with a complicated kernel, they reduced equation $(\ref{ieec})$ to two Fredholm integral equations of the second kind with a simple kernel. Based on these results, Love \cite{love} obtained an infinite series formula for the capacity of a ring, with a certain restriction on the radii of the ring. Our next statement generalizes the known results to the case of higher dimensions and general Riesz $s$-potentials.
\begin{theorem}\label{theo5}
Let $s=(d-3)+2\lambda$, with $0<\lambda<1$. Suppose that an external field $Q$ is invariant with respect to rotations about the polar axis, that is $Q(x)=Q(r)$, where $x=(0,r\overbar{x})\in\D$, $\overbar{x}\in\S^{d-2}$, $0\leq r \leq 1$, and is such that $Q\in C^2(\D)$. Assume that $\supp\mu_Q$ is a ring $\mathscr{R}(a,b)$. Let
\begin{equation}\label{auxfF}
F(r) = \frac{\Gamma((d+2\lambda-3)/2)\Gamma(3-2\lambda)}{2\, \Gamma((d-2\lambda+3)/2)}\, \frac{d}{dr} \int_a^r \frac{Q(t)\,t^{d-2}\,dt}{(r^2-t^2)^{1-\lambda}}, \quad a \leq r \leq b,
\end{equation}
and further put
\begin{empheq}{align}\label{kernel-K}
K(u,r)  =\frac{a^{d-2\lambda+1}}{u^2-r^2} \, \bigg\{ & \frac{1}{r^2}\,  {}_2 F_1\left(1,\frac{d+2\lambda-3}{2};\frac{d-2\lambda+3}{2};\left(\frac{a}{r}\right)^2\right) 
\\\nonumber
& -  \frac{1}{u^2}\,  {}_2 F_1\left(1,\frac{d+2\lambda-3}{2};\frac{d-2\lambda+3}{2};\left(\frac{a}{u}\right)^2\right) \bigg\}.
\end{empheq}
Let $f$ be the density of the extremal measure $\mu_Q$, that is $d\mu_Q(x) = f(r)\, r^{d-2}\, dr\, d\sigma_{d-1}(\overbar{x})$, $x=(0,r\overbar{x})\in\mathscr{R}(a,b)$, $\overbar{x}\in\S^{d-2}$, $a\leq r \leq b$. We set
\begin{equation}\label{density-transform}
G(r) = \int_r^b \frac {f(t)\, t\, dt} {(t^2-r^2)^\lambda}, \quad a \leq r \leq b.
\end{equation}
Then the function $G$ can be recovered by solving the Fredholm integral equation of the second kind,
\begin{empheq}{align}\label{integ-eq-1}
G(r) & -  \frac{\Gamma((d+2\lambda-3)/2)\Gamma(3-2\lambda)}{2\, \Gamma((d-2\lambda+3)/2)}\, \int_a^b  G(u)\, K(u,r)\,du \\\nonumber 
& = F_Q \frac{\Gamma((d+2\lambda-3)/2)}{2\pi^{(d-1)/2}\, \Gamma(\lambda)}\, \bigg\{  \frac{d+2\lambda-3}{2} \, r^{d+2\lambda-4}\, \Beta \left(1-\left(\frac{a}{r}\right)^2;\lambda, \frac{d-1}{2}\right) + \frac {a^{d-1}}{r}\, (r^2-a^2)^{\lambda-1} \bigg\} \\\nonumber 
& - F(r), \quad a \leq r \leq b.
\end{empheq}
The constant $F_Q$ is uniquely determined by the relation
\begin{equation}\label{const-F-Q-ring}
\int_a^b f(t)\, t^{d-2}\, dt = \frac{\Gamma((d-1)/2)}{2\pi^{(d-1)/2}}.
\end{equation}
\end{theorem}

\noindent We note that equation $(\ref{integ-eq-1})$ is a Fredholm integral equation of the second kind. The kernel $K$, given by $(\ref{kernel-K})$, is symmetric, that is $K(u,r)=K(r,u)$, which directly follows from its expression $(\ref{kernel-K})$. A closed-form exact solution of integral equations of type $(\ref{integ-eq-1})$ is presently unknown. However, various numerical methods and methods for finding approximate solutions to Fredholm integral equations of the second kind are well-established. These results, among many other facts on virtually all known types of integral equations can be found in book \cite{polman}.


\section{Proofs}


\noindent{\bf Proof of Theorem \ref{theo0}.} From the rotational invariance of $Q$ it follows that the support of the extremal measure $\mu_Q$ is also rotationally invariant. Therefore, there exists a compact set $K\subset [0,1]$ and a non-negative real-valued function $f\in L^1(K)$ such that 
\begin{align*}
& d\mu_Q(x) = f(r)\, r^{d-2}\, dr\, d\sigma_{d-1}(\overbar{x}), \quad x=(0,r\overbar{x}),\quad \overbar{x}\in\S^{d-2},\\
& \supp\mu_Q=\{(0,r\overbar{x})\in\D: r\in K, \overbar{x}\in\S^{d-2} \}.
\end{align*}
We will first prove that the set $K$ is connected. We will follow the argument given in \cite{bds1}. Assume to the contrary that $K$ is not connected. Then there is an interval $[r_1,r_2]\subset [0,1]$ such that $K\cap [r_1,r_2]=\{r_1,r_2\}$. We further denote $K_+:=K\cap[r_2,1]$ and $K_-:=K\cap[0,r_1]$. Then for
\begin{align*}
& x=(0,r\overbar{x}), \quad r\in(r_1,r_2), \quad \overbar{x}\in\S^{d-2},\\
& y=(0,\rho\overbar{y}), \quad \rho\in K_-\cup K_+, \quad \overbar{y}\in\S^{d-2},
\end{align*}
the $s$-potential of $\mu_Q$ can be written as
\begin{empheq}{align*}\label{newtpot}
U_s^{\mu_Q}(x) & = \int\frac{1}{|x-y|^s}\,d\mu(y)\\ \nonumber
			& = \int_K f(\rho)\, \rho^{d-2}\, d\rho \, \int_{\S^{d-2}} \frac { d \sigma_{d-1}(\overbar{y}) } {(r^2 + \rho^2 - 2r\rho \langle \overbar{x},\overbar{y}\rangle)^{s/2} }\\\nonumber
			& = \frac {2\pi^{(d-2)/2}} {\Gamma(d/2-1)}\, \int_K f(\rho)\, \rho^{d-2}\, d\rho \, \int_0^\pi \frac {\sin^{d-3}\xi \,d\xi} {(r^2 + \rho^2 - 2r\rho \cos\xi)^{s/2}} \\\nonumber
			& = \int_{K_-}  f(\rho) \, k(r,\rho)\, d\rho +  \int_{K_+}  f(\rho) \, k(r,\rho)\, d\rho,
\end{empheq}
where $k(r,\rho)$ is given by
\begin{equation}\label{kern-supp}
k(r,\rho) = \frac {2\pi^{(d-2)/2}} {\Gamma(d/2-1)}\, \int_0^\pi \frac {\sin^{d-3}\xi \,d\xi} {(r^2 + \rho^2 - 2r\rho \cos\xi)^{s/2}}.
\end{equation}
The result  \#3.665 of \cite{Grad} states that for $\operatorname{Re}(\nu)>0$ and $|x|<1$,
\begin{equation}\label{trig-int-rep}
 \int_0^\pi \frac {\sin^{2\nu-1}\xi \,d\xi} {(1 + x^2 + 2x\cos\xi)^s} = \frac {\Gamma(\nu)\sqrt{\pi}}{\Gamma(\nu+1/2)}\, {}_2 F_1(s,s-\nu+1/2;\nu+1/2;x^2).
\end{equation}
Using $(\ref{trig-int-rep})$, for the case $r>\rho$ ($\rho\in K_-$) we can further transform $(\ref{kern-supp})$ as follows,
\begin{align*}
k(r,\rho) & = \frac {2\pi^{(d-1)/2}} {\Gamma((d-1)/2)}\, \frac{\rho^{d-2}} {r^s}\,  {}_2 F_1\left(\frac{s}{2},\lambda;\frac{d-1}{2};\left(\frac{\rho}{r}\right)^2\right) \\
		& = \frac {2\pi^{(d-1)/2}\,\rho^{d-2}} {\Gamma((d-1)/2)}\, \sum_{n=0}^\infty \frac {(s/2)_n\, (\lambda)_n\, \rho^{2n}} {((d-1)/2)_n\, n!} \frac{1} {r^{2n+s}}.
\end{align*}
Hence, for $r>\rho$ we have
\begin{equation}\label{k-series}
k(r,\rho) = \frac {2\pi^{(d-1)/2}\,\rho^{d-2}} {\Gamma((d-1)/2)}\, \sum_{n=0}^\infty \frac {(s/2)_n\, (\lambda)_n\, \rho^{2n}} {((d-1)/2)_n\, n!}\, \frac{1} {r^{2n+s}}.
\end{equation}
It is clear that the functions $r^{-(2n+s)}$, $n=0,1,2,\ldots$ are strictly convex for  $r\in (0,1)$. Therefore, taking into account the positivity of all coefficients of the series in the right hand side of $(\ref{k-series})$, and the fact that this series is uniformly convergent in $r$ on compact subsets of $(r_1,r_2)$, the convexity of the right hand side of $(\ref{k-series})$ follows by differentiation with respect to $r$.

Exactly the same approach with $\rho>r$ ($\rho\in K_+$) leads to the following series representation for $k(r,\rho)$,
\begin{equation}\label{k+series}
k(r,\rho) = \frac {2\pi^{(d-1)/2}} {\Gamma((d-1)/2)\,\rho^{2\lambda-1}}\, \sum_{n=0}^\infty \frac {(s/2)_n\, (\lambda)_n} {((d-1)/2)_n\, n!\, \rho^{2n}}\, r^{2n}.
\end{equation}
Obviously the functions $r^{2n}$, $n=0,1,2\ldots$ are convex on $(0,1)$. Hence, we similarly derive the convexity of the right hand side of $(\ref{k+series})$.

From $(\ref{k-series})$ and $(\ref{k+series})$ we infer that the function $k$ is a strictly convex function of $r$ on $(r_1,r_2)$, for any fixed $\rho\in K_- \cup K_+$. Using the convexity of $Q(r)$, we deduce that $U_s^{\mu_Q}(r)+Q(r)$ is  a strictly convex function on $(r_1,r_2)$. Furthermore,  by $(\ref{var2})$, for the weighted potential $U_s^{\mu_Q}(r)+Q(r)$ we have $U_s^{\mu_Q}(r_1)+Q(r_1)=F_Q=U_s^{\mu_Q}(r_2)+Q(r_2)$. Then the strict convexity of $U_s^{\mu_Q}(r)+Q(r)$ implies that $U_s^{\mu_Q}(r)+Q(r)<F_Q$, for $r_1<r<r_2$. But this is an obvious contradiction with inequality $(\ref{var1})$, which is valid for all $0\leq r\leq 1$.

We now prove the second part of the statement of Theorem \ref{theo0}. Assume that $Q(r)$, in addition to being convex, is also an increasing function. Suppose that $a>0$. In this case the kernel $k(r,\rho)$ is calculated according to $(\ref{k+series})$, which shows that $U_s^{\mu_Q}(r)$ is an increasing function on $(0,a]$. This implies that the weighted potential $U_s^{\mu_Q}(r)+Q(r)$ is a strictly increasing function on $(0,a]$. Therefore for any $r'\in(0,a)$ we have $U_s^{\mu_Q}(a)+Q(a)>U_s^{\mu_Q}(r')+Q(r')$. On the other hand, since $a\in S_Q$, from $(\ref{var2})$ it follows that $U_s^{\mu_Q}(a)+Q(a)=F_Q$. We thus find that $U_s^{\mu_Q}(r')+Q(r')<F_Q$, which clearly violates inequality $(\ref{var1})$.

The proof of the remaining part of the statement of Theorem \ref{theo0} follows the same logic as in the last paragraph. Indeed, assume that $Q(r)$ besides being convex, is also a decreasing function of $r$. Suppose that $b<1$. In this case kernel $k(r,\rho)$ is calculated according to $(\ref{k-series})$, which shows that $U_s^{\mu_Q}(r)$ is a decreasing function on $[b,1)$. Hence the weighted potential $U_s^{\mu_Q}(r)+Q(r)$ is a strictly decreasing function on $[b,1)$. Thus for any $r'\in(b,1)$ we have $U_s^{\mu_Q}(b)+Q(b)>U_s^{\mu_Q}(r')+Q(r')$. Observing that $b\in S_Q$, from $(\ref{var2})$ it follows that $U_s^{\mu_Q}(b)+Q(b)=F_Q$. We thus see that  $U_s^{\mu_Q}(r')+Q(r')<F_Q$, which again violates inequality $(\ref{var1})$.

\qed


\noindent{\bf Proof of Theorem \ref{theo1}.} Assume that the support of the extremal measure $\mu_Q$ is the disk $\D_R$, that is $S_Q = \D_R$. Then there exists a non-negative real-valued function $f\in L^1([0,R])$ such that for $x=(0,r\overbar{x})\in\D_R$, $\overbar{x}\in\S^{d-2}$, $0\leq r \leq R$, 
\[
d\mu_Q(x) = f(r)\, r^{d-2}\, dr\, d\sigma_{d-1}(\overbar{x}).
\]
\noindent Let $x$ and $y$ be two points in $\D$ with $|x|=r$ and $|y|=\rho$. Note that $\overbar{x}:=x/r\in\S^{d-2}$, and similarly $\overbar{y}:=y/\rho\in\S^{d-2}$. For the distance $|x-y|$ we then obtain
\begin{align*}
 |x-y|^2 & = |x|^2 + |y|^2 - 2 \langle x,y \rangle \\
 		& = r^2 + \rho^2 - 2r\rho \langle \overbar{x},\overbar{y}\rangle.
\end{align*}
We immediately notice that the rotational invariance of an external field $Q$ is passed on to the Riesz $s$-potential $U_s^{\mu_Q}$, thanks to the uniqueness of the extremal measure. Therefore for $x\in\supp\mu_Q$ we have $U_s^{\mu_Q}(x)=U_s^{\mu_Q}(|x|)$. We will be using this fact from now on without mentioning it explicitly on each separate occasion.

The Riesz $s$-potential $U_s^{\mu_Q}(x)$, with $x=(0,r\overbar{x})\in\D_R$, $\overbar{x}\in\S^{d-2}$, $0\leq r \leq R$, can be written as
\begin{empheq}{align}\label{pinspc}
U_s^{\mu_Q}(x) & = \int_{\D_R} \frac {1} {|x-y|^s}\, d\mu_Q(y) \\ \nonumber
			 & = \int_0^R  f(\rho) \, \rho^{d-2} d\rho \, \int_{\S^{d-2}} \frac { d \sigma_{d-1}(\overbar{y}) } {(r^2 + \rho^2 - 2r\rho \langle \overbar{x},\overbar{y}\rangle)^{s/2} }.
\end{empheq}
We will need the following proposition, which is a special case of the Funk-Hecke theorem \cite[p. 247]{be2}.
\begin{proposition}\label{funkhecke}
If $f$ is integrable on $[-1,1]$ with respect to the weight $(1-t^2)^{(d-3)/2}$, and $y$ is an arbitrary fixed point on the sphere $S^{d-1}$, then
\begin{equation}\label{fhformula}
\int_{\S^{d-1}} f (\langle x,y \rangle ) \, d\sigma_d(x) = \frac {2 \pi^{(d-1)/2}} {\Gamma ((d-1)/2)} \, \int_{-1}^1 f(t) \, (1-t^2)^{(d-3)/2} \, dt.
\end{equation}
\end{proposition}

\noindent Applying Proposition \ref{funkhecke} to the inner integral on the right hand side in $(\ref{pinspc})$, we derive
\begin{equation*}
 \int_{\S^{d-2}} \frac { d \sigma_{d-1}(\overbar{y}) } {(r^2 + \rho^2 - 2r\rho \langle \overbar{x},\overbar{y}\rangle)^{s/2}} =  \frac {2 \pi^{(d-2)/2}} {\Gamma ((d-2)/2)} \, \int_0^\pi \frac {\sin^{d-3}\xi \,d\xi} {(r^2 + \rho^2 - 2r\rho \cos\xi)^{s/2}}.
\end{equation*}
Hence the potential $U_s^{\mu_Q}$ in $(\ref{pinspc})$ assumes the form
\[
U_s^{\mu_Q}(r)  =   \frac {2 \pi^{(d-2)/2}} {\Gamma ((d-2)/2)} \,  \int_0^R  f(\rho) \, \rho^{d-2} d\rho \, \int_0^\pi \frac {\sin^{d-3}\xi \,d\xi} {(r^2 + \rho^2 - 2r\rho \cos\xi)^{s/2}}, \quad 0 \leq r \leq R.
\]
The integral equation $(\ref{ieec})$ can now be written as
\begin{equation}\label{ie1}
 \int_0^R  f(\rho) \, \rho^{d-2} d\rho \, \int_0^\pi \frac {\sin^{d-3}\xi \,d\xi} {(r^2 + \rho^2 - 2r\rho \cos\xi)^{s/2}}  =   \frac {\Gamma ((d-2)/2)} {2 \pi^{(d-2)/2}}  \,(F_Q - Q(r)), \quad 0 \leq r \leq R.
\end{equation}
Our next step is to further transform the inner integral on the left hand side of the integral equation $(\ref{ie1})$. This is achieved via the following fact.
\begin{lemma}\label{integ-rep}
If $a$ and $b$ are positive numbers, $a\neq b$, $q \geq 0$ and $0<\lambda<1$, then
\begin{equation}\label{integ-rep-form}
\int_0^{\pi} \frac {\sin^{2q}\xi \,d\xi} {(a^2+b^2-2ab\cos\xi)^{q+\lambda}} = \frac{2}{a^{2q}\, b^{2q}}\, \frac{\sin(\lambda\pi)\, \Gamma(\lambda)\, \Gamma(q+1/2)} {\sqrt{\pi}\, \Gamma(q+\lambda)}\, \int_0^{\min{(a,b)}} \frac {t^{2(q+\lambda)-1}\,dt} {(a^2-t^2)^\lambda\, (b^2-t^2)^\lambda}.
\end{equation}
\end{lemma}
\noindent Note that Lemma \ref{integ-rep} in the case when $q=0$ and $\lambda=1/2$ was obtained by Copson \cite{cop}. Also, Lemma \ref{integ-rep} when $q\geq0$ and $\lambda=1/2$, is implicitly mentioned in \cite{shail}, although with an incorrect numerical coefficient. The correct version of the latter, along with its proof, is given in \cite[p. 8]{bilo2}.

\begin{proof}
The proof is based on the following identity, obtained by Kahane \cite{kah}.
\begin{proposition}\label{kahlemma}
Let $a$ and $b$ be positive numbers such that $a\neq b$, $\lambda\in(0,1)$, $q\in\mathbb C$ with $\operatorname{Re}(q)\geq0$, and $u$ a real number with $|u|\leq 1$. Then
\begin{empheq}{align}\label{kahformula}
 \frac{(ab)^q} {(a^2+b^2-2abu)^{\lambda+q}} & = \frac {\Gamma(\lambda) \Gamma(q+1)} {\Gamma(\lambda+q)} \frac {2\sin(\lambda \pi)} {\pi} \times   \\ \nonumber
& \int_0^{\min{(a,b)}} \frac {1-(t^2/ab)^2} {(1+(t^2/ab)^2-2(t^2/ab)u)^{q+1}} \left( \frac {t^2} {ab} \right)^q \frac {t^{2\lambda-1} \,dt} {(a^2-t^2)^{\lambda} (b^2-t^2)^{\lambda}}.
\end{empheq}
\end{proposition}
\noindent Taking $q$ to be a non-negative real number in Proposition \ref{kahlemma}, and applying Fubini's theorem, we rewrite the left hand side of $(\ref{integ-rep-form})$ as
\begin{empheq}{align}\label{kern1}
\int_0^{\pi} \frac {\sin^{2q}\xi \,d\xi} {(a^2+b^2-2ab\cos\xi)^{q+\lambda}}  =  \frac{1}{a^{2q}\, b^{2q}} \frac {2\sin(\lambda\pi)\,\Gamma(\lambda)\,\Gamma(q+1)} {\pi\, \Gamma(q+\lambda)} & \int_0^{\min{(a,b)}}  \frac {(1-(t^2/ab)^2)\, t^{2(q+\lambda)-1}\,dt} { (a^2-t^2)^\lambda\, (b^2-t^2)^\lambda} \times \\ \nonumber
& \int_0^\pi \frac{\sin^{2q}\xi \, d\xi} {(1+(t^2/ab)^2-2(t^2/ab)\cos\xi)^{q+1}}.
\end{empheq}
We will show that
\begin{equation}\label{ltint}
\int_0^\pi \frac{\sin^{2q}\xi \, d\xi} {(1+(t^2/ab)^2-2(t^2/ab)\cos\xi)^{q+1}} = \frac {1} {1-(t^2/ab)^2} \frac {\sqrt{\pi}\, \Gamma(q+1/2)} {\Gamma(q+1)}.
\end{equation}
Indeed, the integral of a type appearing on the left hand side of $(\ref{ltint})$ was previously considered in \cite[p. 400]{land}. It was shown that
\begin{equation}\label{lankint}
\int_0^\pi \frac {\sin^{p-2}\xi \, d\xi} {(1+\rho^2-2\rho\cos\xi)^{p/2}} = \frac {1} {\rho^{p-2}(\rho^2-1)} \int_0^\pi \sin^{p-2}\xi \,d\xi,
\end{equation}
where $\rho\geq 1$, and $p\geq 3$ was assumed to be an integer. A careful analysis of the evaluation of integral $(\ref{lankint})$ in \cite[p. 400]{land} shows that, in fact, $(\ref{lankint})$ holds true
for any $p\geq2$.
We hence transform the left hand side of $(\ref{ltint})$ as follows,
\begin{align*}
\int_0^\pi \frac{\sin^{2q}\xi \, d\xi} {(1+(t^2/ab)^2-2(t^2/ab)\cos\xi)^{q+1}} & = \frac {1} {1-(t^2/ab)^2} \int_0^\pi \sin^{2q}\xi \, d\xi \\
			& = \frac {1} {1-(t^2/ab)^2} 2^{2q} \, \Beta(q+1/2,q+1/2) \\
			& = \frac {1} {1-(t^2/ab)^2}  \frac {\sqrt{\pi}\, \Gamma(q+1/2)} {\Gamma(q+1)},
\end{align*}
which is the right hand side of $(\ref{ltint})$. Substituting $(\ref{ltint})$ into $(\ref{kern1})$, we obtain the desired representation $(\ref{integ-rep-form})$.
\end{proof}

Setting $q=(d-3)/2$ and observing that $s=2q+2\lambda$, by letting $a=r$ and $b=\rho$ in Lemma \ref{integ-rep-form}, the inner integral on the left hand side of integral equation $(\ref{ie1})$ can be written in the following form,
\begin{equation}\label{int-rep-kern}
\int_0^\pi \frac {\sin^{d-3}\xi \,d\xi} {(r^2 + \rho^2 - 2r\rho \cos\xi)^{s/2}} = \frac{1}{r^{d-3}\, \rho^{d-3}} \, \frac {2\sin(\lambda\pi)\, \Gamma(\lambda)\, \Gamma(d/2-1)}{\sqrt{\pi}\, \Gamma((d-3)/2+\lambda)} \, \int_0^{\min{(r,\rho)}} \frac {t^{d+2\lambda-4}\,dt} {(r^2-t^2)^\lambda \, (\rho^2-t^2)^\lambda}.
\end{equation}
Using $(\ref{int-rep-kern})$, we recast integral equation $(\ref{ie1})$ as
\begin{equation}\label{ie2}
\int_0^R f(\rho)\, \rho \, d\rho \int_0^{\min{(r,\rho)}} \frac {t^{d+2\lambda-4}\,dt} {(r^2-t^2)^\lambda \, (\rho^2-t^2)^\lambda} = \frac {\Gamma((d-3)/2+\lambda)}{4\sin(\lambda\pi)\, \pi^{(d-3)/2}\, \Gamma(\lambda)}\, r^{d-3}\,(F_Q-Q(r)), \quad 0\leq r\leq R.
\end{equation}
We now work with the integral on the left hand side of $(\ref{ie2})$. Splitting the range of integration, and changing the order of integration in the first integral, we derive
\begin{align*}
\int_0^R f(\rho)\, \rho \, d\rho \int_0^{\min{(r,\rho)}} \frac {t^{d+2\lambda-4}\,dt} {(r^2-t^2)^\lambda \, (\rho^2-t^2)^\lambda}  = & \int_0^r f(\rho)\, \rho \, d\rho \int_0^{\rho} \frac {t^{d+2\lambda-4}\,dt} {(r^2-t^2)^\lambda \, (\rho^2-t^2)^\lambda}  \\
+ & \int_r^R f(\rho)\, \rho \, d\rho \int_0^{r} \frac {t^{d+2\lambda-4}\,dt} {(r^2-t^2)^\lambda \, (\rho^2-t^2)^\lambda} \\
=  & \int_0^r \frac {t^{d+2\lambda-4}\,dt} {(r^2-t^2)^\lambda}  \int_t^r \frac{f(\rho)\, \rho \, d\rho} {(\rho^2-t^2)^\lambda} \\
+ & \int_0^r \frac {t^{d+2\lambda-4}\,dt} {(r^2-t^2)^\lambda}  \int_r^R \frac{f(\rho)\, \rho \, d\rho} {(\rho^2-t^2)^\lambda} \\
= & \int_0^r \frac {t^{d+2\lambda-4}\,dt} {(r^2-t^2)^\lambda}  \int_t^R \frac{f(\rho)\, \rho \, d\rho} {(\rho^2-t^2)^\lambda}.
\end{align*}
We can thus re-write integral equation $(\ref{ie2})$ as
\begin{equation}\label{ie3}
\int_0^r \frac {t^{d+2\lambda-4}\,dt} {(r^2-t^2)^\lambda}  \int_t^R \frac{f(\rho)\, \rho \, d\rho} {(\rho^2-t^2)^\lambda} = \frac {\Gamma((d-3)/2+\lambda)}{4\sin(\lambda\pi)\, \pi^{(d-3)/2}\, \Gamma(\lambda)}\, r^{d-3}\,(F_Q-Q(r)), \quad 0\leq r\leq R.
\end{equation}
Let
\begin{equation}\label{func-s}
S(t) = \int_t^R \frac{f(\rho)\, \rho \, d\rho} {(\rho^2-t^2)^\lambda}.
\end{equation}
Then  $(\ref{ie3})$ reads
\begin{equation}\label{ie4}
\int_0^r \frac {S(t)\, t^{d+2\lambda-4}\,dt} {(r^2-t^2)^\lambda} = \frac {\Gamma((d-3)/2+\lambda)}{4\sin(\lambda\pi)\, \pi^{(d-3)/2}\, \Gamma(\lambda)}\, r^{d-3}\,(F_Q-Q(r)), \quad 0\leq r\leq R.
\end{equation}
Integral equation $(\ref{ie4})$ is an Abel-type integral equation with respect to $S(t)\, t^{d+2\lambda-4}$. As $Q\in C^2([0,1])$, applying \cite[\# 44, p. 122]{polman}, we solve this equation and find
\begin{equation}\label{func-s-sol}
S(r) = \frac {\Gamma((d-3)/2+\lambda)} {2\, \pi^{(d-1)/2}\, \Gamma(\lambda)}\, \frac {1}{r^{d+2\lambda-4}}\, \frac{d}{dr} \int_0^r \frac {(F_Q-Q(t))\, t^{d-2}\, dt} {(r^2-t^2)^{1-\lambda}}.
\end{equation}
Now observe that $(\ref{func-s})$ is also an Abel-type integral equation with respect to $f(\rho)\, \rho$. Solving it in a similar fashion, we derive
\begin{equation}\label{ie-sol}
f(t) = - \frac {2\sin(\lambda\pi)}{\pi}\, \frac{1}{t}\, \frac{d}{dt} \int_t^R \frac {S(\rho)\, \rho\, d\rho} {(\rho^2-t^2)^{1-\lambda}}, \quad 0\leq t \leq R.
\end{equation}
Let
\begin{equation*}
F(t) = \frac {\sin(\lambda\pi)\, \Gamma((d-3)/2+\lambda)}{\pi^{(d+1)/2}\, \Gamma(\lambda)}\, \frac{1}{t}\, \frac{d}{dt} \int_t^R \frac {g(r)\, r\, dr} {(r^2-t^2)^{1-\lambda}}, \quad 0\leq t \leq R,
\end{equation*}
where
\begin{equation*}
g(r) = \frac {1}{r^{d+2\lambda-4}}\, \frac{d}{dr} \int_0^r \frac {Q(u)\, u^{d-2}\, du} {(r^2-u^2)^{1-\lambda}}, \quad 0\leq r \leq R.
\end{equation*}
Then expression $(\ref{ie-sol})$ can be written as
\begin{align*}
f(t) =  - & F_Q\, \frac {\sin(\lambda\pi)\, \Gamma((d-3)/2+\lambda)}{\pi^{(d+1)/2}\, \Gamma(\lambda)}\, \frac{1}{t}\, \frac{d}{dt} \int_t^R \bigg\{ \frac {1}{r^{d+2\lambda-4}}\, \frac{d}{dr} \int_0^r \frac {u^{d-2}\, du} {(r^2-u^2)^{1-\lambda}} \bigg\} \frac {r\, dr} {(r^2-t^2)^{1-\lambda}} \\
+ & F(t), \quad 0\leq t \leq R.
\end{align*}
Performing rather straightforward integrations and differentiations appearing on the right hand side of the last expression, we deduce
\begin{equation}\label{extr-density}
f(t) = F_Q\, \frac {\sin(\lambda\pi)\, \Gamma((d-1)/2)} {\pi^{(d+1)/2}}\, (R^2-t^2)^{\lambda-1} + F(t), \quad 0\leq t \leq R.
\end{equation}
We complete the proof by evaluating the Robin constant $F_Q$. Recall that $\mu_Q$ is a probability measure, so that it has mass one. Therefore,
\begin{align*}
1=\int d\mu_Q & = \int_0^R f(t)\, t^{d-2}\, dt\, \int_{\S^{d-2}} d\sigma_{d-1} \\
				& = \omega_{d-1}\, \int_0^R f(t)\, t^{d-2}\, dt.
\end{align*}
We therefore obtain
\begin{empheq}{align*}\label{constfq}
 \frac{\Gamma((d-1)/2)}{2\pi^{(d-1)/2}}  = \frac{1}{\omega_{d-1}}  &  =  \int_0^R f(t)\, t^{d-2}\, dt \\ \nonumber
	& =F_Q\, \frac{\sin(\lambda\pi)\, \Gamma((d-1)/2)}{\pi^{(d+1)/2}}\, \int_0^R(R^2-t^2)^{\lambda-1}\, t^{d-2}\, dt \\ \nonumber
	& + \int_0^R F(t)\, t^{d-2}\, dt.
\end{empheq}
It is an elementary calculation to see that
\begin{equation}\label{lambda-int-1}
\int_0^R(R^2-t^2)^{\lambda-1}\, t^{d-2}\, dt = \frac {\Gamma((d-1)/2)\, \Gamma(\lambda)} {2\, \Gamma((d-1)/2+\lambda)}\, R^{d+2\lambda-3}.
\end{equation}
Combining the last two expressions, we eventually find
\begin{equation}\label{rob-const-F}
F_Q = \frac {2\, \pi^{(d+1)/2}\, \Gamma((d-1)/2+\lambda)} {\sin(\lambda\pi)\, (\Gamma((d-1)/2))^2\, \Gamma(\lambda)}\, \frac{1}{R^{d+2\lambda-3}} \, \bigg\{ \frac {\Gamma((d-1)/2)} {2\, \pi^{(d-1)/2}} - \int_0^R F(t)\, t^{d-2}\, dt \bigg\}.
\end{equation}
Inserting expression $(\ref{rob-const-F})$ into the formula for the extremal density $(\ref{extr-density})$, we derive the following simple formula,
\begin{equation*}
f(t) = C_Q\, (R^2-t^2)^{\lambda-1} + F(t), \quad 0\leq t \leq R,
\end{equation*}
where the constant $C_Q$ is given by
\begin{equation*}
C_Q = \frac {2\, \Gamma((d-1)/2+\lambda)} {\Gamma(\lambda)\, \Gamma((d-1)/2)}\, \frac{1}{R^{d+2\lambda-3}} \, \bigg\{ \frac {\Gamma((d-1)/2)} {2\, \pi^{(d-1)/2}} - \int_0^R F(t)\, t^{d-2}\, dt \bigg\}.
\end{equation*}
\qed


\noindent{\bf Proof of Theorem \ref{equil-meas-cap-disk}.} The expression for the equilibrium measure on the disk $\D_R$ when there is no external field present, follows from Theorem \ref{theo1} upon setting $Q=0$. To obtain expression for the capacity of the disk $\D_R$, we first notice that $F_Q=W_s(\D_R)$ when $Q=0$. From formula 
$(\ref{rob-const-F})$ we find that when $Q=0$,
\begin{equation*}
F_Q = \frac {\pi\, \Gamma((d+2\lambda-1)/2)} {\sin(\lambda\pi)\,  \Gamma(\lambda)\, (\Gamma((d-1)/2))^2}\, \frac{1}{R^{d+2\lambda-3}}.
\end{equation*}
Recalling that $\cp_s(\D_R)=1/W_s(\D_R)$, we obtain desired expression $(\ref{cap-disk})$.

\qed


\noindent{\bf Proof of Proposition \ref{ms-func-min-prop}.} The proof below follows the same line of argument as in \cite{bds1}. Let $E$ be any compact subset of $\D$ with positive Riesz $s$-capacity. For the range of the Riesz $s$-parameter satisfying $d-3<s<d-1$, the potential of the equilibrium measure $\mu_E$ (with no external field) satisfies the following inequalities \cite[p. 136]{land},
\begin{empheq}{align}\label{viu1}
& U_s^{\mu_E}(x) = W_s(E), \quad \text{q. e. on}\ \ E,\\\label{viu2}
& U_s^{\mu_E} (x) \leq W_s(E), \quad \text{on} \ \ \D.
\end{empheq}
We first observe that variational inequalities $(\ref{var1})$--$(\ref{var2})$ imply
\begin{empheq}{align*}
 \cF_s(S_Q) & = W_s(S_Q) + \int Q(x)\,d\mu_Q(x)  \\
 		& = I_s(\mu_Q) +\int Q(x)\,d\mu_Q(x) \\
		& = F_Q.
\end{empheq}
We now show that for any compact set $E\subset\D$ with positive Riesz $s$-capacity we have $\cF_s(E)\geq \cF_s(S_Q)$. Indeed, integrating inequality $(\ref{var1})$ with respect to $\mu_E$, we obtain
\begin{equation}\label{auxineq}
\int U_s^{\mu_Q}(x)\,d\mu_{E}(x)  + \int Q(x)\,d\mu_{E}(x)  \geq F_Q,
\end{equation}
where the inequality holds $\mu_E$-a.e. as $\mu_E$ has finite Riesz $s$-energy. With $(\ref{auxineq})$ and $(\ref{viu1})$--$(\ref{viu2})$ in mind, we write down the following chain of inequalities,
\begin{empheq}{align*}
W_s(E) & = \int W_s(E)\,d\mu_Q(x)\\
	& \geq \int U_s^{\mu_E}(x) \,d\mu_Q(x) \\
	& = \int \left( \int \frac {1} {|x-y|^s} \,d\mu_E(y) \right)\,d\mu_Q(x)\\
	& = \int \left( \int \frac {1} {|x-y|^s} \,d\mu_Q(x) \right)\,d\mu_E(y)\\
	& = \int U_s^{\mu_Q}(y) \,d\mu_E(y)\\
	& = \int U_s^{\mu_Q}(x) \,d\mu_E(x)\\
	& \geq F_Q - \int Q(x)\,d\mu_{E}(x).
\end{empheq}
We now see that
\begin{equation*}
\cF_s(E) = W_s(E) + \int Q(x)\,d\mu_{E}(x) \geq F_Q = \cF_s(S_Q), 
\end{equation*}
so that $\cF_s(E) \geq \cF_s(S_Q)$, as claimed.

\qed


\noindent{\bf Proof of Theorem \ref{critical-radius-theo}.} If $E=\D_R$, taking into account that $W_s(\D_R)=1/\cp_s(\D_R)$, and inserting $(\ref{cap-disk})$ and $(\ref{equil-meas-disk})$ into $(\ref{ms-func})$, we find that $\cF_s$-functional is given by
\begin{equation}\label{ms-func-disk-R}
 \cF_s(\D_R) =  \frac   {\pi\, \Gamma((d+2\lambda-1)/2)} {\sin(\lambda\pi)\, \Gamma(\lambda)\, \Gamma((d-1)/2)}\, \frac{1}{R^{d+2\lambda-3}} \bigg\{  1+ \frac{2\sin(\lambda\pi)}{\pi}\, \int_0^R Q(r)\, (R^2-r^2)^{\lambda-1}\, r^{d-2}\, dr \bigg\}.
\end{equation}
Using the substitution $R-r=Ru$, we transform the integral on the right hand side of $(\ref{ms-func-disk-R})$ as follows,
\begin{equation*}
\int_0^R Q(r) (R^2-r^2)^{\lambda-1} \, r^{d-2}\, dr = 2^{\lambda-1}\, R^{d+2\lambda-3}\, \int_0^1 Q((1-u)R)\,u^{\lambda-1}\, (1-u)^{d-2}\, (1-u/2)^{\lambda-1}\,du.
\end{equation*}
This allows us to write expression $(\ref{ms-func-disk-R})$ as
\begin{equation*}
\cF_s(\D_R) =  c(d,\lambda)\, \bigg\{\frac{1}{R^{d+2\lambda-3}} + \frac{2^\lambda \sin(\lambda\pi)}{\pi}\,\int_0^1 Q((1-u)R)\,u^{\lambda-1}\, (1-u)^{d-2}\, (1-u/2)^{\lambda-1}\,du \bigg\}, 
\end{equation*}
where for brevity we set
\begin{equation*}
c(d,\lambda) := \frac {\pi \Gamma((d+2\lambda-1)/2)} {\sin(\lambda\pi) \Gamma(\lambda) \Gamma((d-1)/2)}.
\end{equation*}
Differentiating the last expression with respect to $R$, we derive
\begin{empheq}{align}\label{ms-func-first-derivat}
\cF_s'(\D_R) & =  c(d,\lambda)\,  \bigg\{ - \frac{d+2\lambda-3}{R^{d+2\lambda-2}} + \frac{2^\lambda \sin(\lambda\pi)}{\pi}\, \int_0^1 Q'((1-u)R)\,u^{\lambda-1}\, (1-u)^{d-1}\, (1-u/2)^{\lambda-1}\,du \bigg\} \\ \nonumber
  & = c(d,\lambda)\,  \bigg\{ - \frac{d+2\lambda-3}{R^{d+2\lambda-2}} + \frac{2\sin(\lambda\pi)}{\pi}\, \frac{1}{R^{d+2\lambda-2}} \, \int_0^R Q'(r)\, (R^2-r^2)^{\lambda-1}\, r^{d-1}\, dr \bigg\} \\\nonumber
  & = -\frac{2\sin(\lambda\pi)c(d,\lambda)}{\pi R^{d+2\lambda-2}}\, \Delta(R),
\end{empheq}
with
\begin{equation*}
\Delta(R):=  \frac{\pi(d+2\lambda-3)}{2\sin(\lambda\pi)} - \int_0^R Q'(r)\, (R^2-r^2)^{\lambda-1}\, r^{d-1}\, dr,
\end{equation*}
and where the differentiation under the integral sign is justified by invoking the Dominated Convergence Theorem. Since the $\cF_s$-functional is minimized on the support of the extremal measure, we obtain by setting $\cF_s'(\D_R)$ to zero that the radius $R$ must satisfy the equation
\begin{equation}\label{radius-equation}
 \int_0^R Q'(r)\, (R^2-r^2)^{\lambda-1}\, r^{d-1}\, dr =\frac{\pi(d+2\lambda-3)}{2\sin(\lambda\pi)},
 \end{equation}
which by a simple rearrangement can be brought into the required form $(\ref{radius-equation-equation})$.

We now discuss the existence and uniqueness of a solution to equation $(\ref{radius-equation})$. First we note that the left hand side can be written as
\begin{equation*}
w(R):= \int_0^R Q'(r)\, (R^2-r^2)^{\lambda-1}\, r^{d-1}\, dr = 2^{\lambda-1}\, R^{d+2\lambda-2}\, \int_0^1 Q'((1-u)R)\, u^{\lambda-1}\, (1-u)^{d-1}\, (1-u/2)^{\lambda-1}\, du.
\end{equation*}
From the convexity of $Q$ it follows that $Q'(R)\geq 0$ for $0\leq R\leq 1$. Using the Dominated Convergence Theorem we can show that the integral on the right hand side of the last expression is an increasing function of $R$ for $0\leq R\leq 1$. As $w(R)$ is a product of two non-negative increasing functions on $[0,1]$, the function $w(R)$ itself is an increasing function of $R$ for $0\leq R\leq 1$. Moreover, a simple calculation shows that $w(0)=0$. These considerations make it clear that in a situation when $(\ref{radius-equation})$ does not have a solution on the interval $[0,1]$, it must be the case that
\begin{align*}
\Delta(R) & = \frac{\pi(d+2\lambda-3)}{2\sin(\lambda\pi)} -  \int_0^R Q'(r)\, (R^2-r^2)^{\lambda-1}\, r^{d-1}\, dr \\
		& = \frac{\pi(d+2\lambda-3)}{2\sin(\lambda\pi)} - g(R) > 0, \quad 0\leq R \leq 1.
\end{align*}
From $(\ref{ms-func-first-derivat})$ it then follows that $\cF_s'(\D_R)<0$ for $0 \leq R \leq 1$. Hence $\cF_s(\D_R)$ is strictly decreasing on $[0,1]$ and attains its global minimum at $R=1$.

Suppose now that equation $(\ref{radius-equation})$ does have a solution on the interval $[0,1]$, i.e. $\cF_s(\D_R)$ has a critical point on  $[0,1]$. From $(\ref{ms-func-first-derivat})$ we deduce that
\begin{empheq}{align*}
\cF_s''(\D_R) & =  c(d,\lambda)\,  \bigg\{  \frac{(d+2\lambda-3)(d+2\lambda-2)}{R^{d+2\lambda-1}} + \frac{2^\lambda \sin(\lambda\pi)}{\pi}\, \int_0^1 Q''((1-u)R)\,u^{\lambda-1}\, (1-u)^{d}\, (1-u/2)^{\lambda-1}\,du \bigg\} \\ \nonumber
  & = c(d,\lambda)\,  \bigg\{ \frac{(d+2\lambda-3)(d+2\lambda-2)}{R^{d+2\lambda-1}} + \frac{2\sin(\lambda\pi)}{\pi}\, \frac{1}{R^{d+2\lambda-1}} \, \int_0^R Q''(r)\, (R^2-r^2)^{\lambda-1}\, r^{d}\, dr \bigg\} \\
& = \frac{c(d,\lambda) 2\sin(\lambda\pi)}{\pi R^{d+2\lambda-1}} \bigg\{ \frac{\pi(d+2\lambda-3)(d+2\lambda-2)}{2\sin(\lambda\pi)} + \int_0^R Q''(r)\, (R^2-r^2)^{\lambda-1}\, r^{d}\, dr \bigg\}.
\end{empheq}
Recalling that $Q(r)$ is convex on $[0,1]$, from last expression it transpires that $\cF_s''(\D_R)>0$ for all $0\leq R\leq 1$. This means that if $R_*$ is a critical point of $\cF_s(\D_R)$, then it will be its global minimum.

We also observe that $R=0$ cannot satisfy equation $(\ref{radius-equation})$, which is evident from $(\ref{radius-equation})$ itself. We therefore conclude that $\cF_s(\D_R)$ has exactly one global minimum on $(0,1]$, which is either the unique solution $R_*\in(0,1]$ of equation $(\ref{radius-equation})$ if it exists, or $R_*=1$.

\qed



\noindent{\bf Proof of Theorem \ref{mon-ext-field-extrem-meas}.} Inserting $(\ref{ext-field-monom-gen})$ into the integral in equation $(\ref{radius-equation-equation})$, we find, using formula $(\ref{lambda-int-1})$,
\begin{equation*}
\int_0^R Q'(r)\, (R^2-r^2)^{\lambda-1}\, r^{d-2}\, dr = \frac {q\alpha\Gamma((d+\alpha-1)/2)\, \Gamma(\lambda)} {2 \Gamma((d+\alpha+2\lambda-1)/2)}\, R^{d+\alpha+2\lambda-3}.
\end{equation*}
Inserting the above result back to equation $(\ref{radius-equation-equation})$, after simple algebra we derive desired formula $(\ref{mon-ext-field-extrem-meas-disk-rad})$.

We now proceed with evaluating the density of the extremal measure $\mu_Q$, as outlined in $(\ref{func-g})$ and $(\ref{func-F})$. Inserting expression for the external field $(\ref{ext-field-monom-gen})$ into  $(\ref{func-g})$, we arrive at the integral
\begin{equation*}
\int_0^r \frac{Q(u)\,u^{d-2}\,du}{(r^2-u^2)^{1-\lambda}} = \frac{q(\alpha+d-3)}{4\lambda}\, \Beta\left(1+\lambda,\frac{d+\alpha-3}{2}\right)\, r^{d+2\lambda-3+\alpha},
\end{equation*}
which is evaluated using the substitution $tr^2=r^2-u^2$. We therefore easily find that
\begin{equation*}
g(r) = \frac{q(\alpha+d-3)(d+2\lambda+\alpha-3)}{4\lambda}\, \Beta\left(1+\lambda,\frac{d+\alpha-3}{2}\right)\, r^\alpha.
\end{equation*}
Inserting this result into $(\ref{func-F})$, we arrive at another integral
\begin{equation*}
\int_t^{R_\ast} \frac{r^\alpha\,r\,dr}{(r^2-t^2)^{1-\lambda}} = \frac{R_\ast^\alpha}{2\lambda}\, (R_\ast^2-t^2)^\lambda\, {}_2 F_1\left(-\frac{\alpha}{2}, 1;\lambda+1; 1-\left(\frac{t}{R_\ast}\right)^2\right),
\end{equation*}
which is handled by substituting $(R_\ast^2-t^2)u=r^2-t^2$ and then recalling the integral representation \cite[\#15.3.1, p. 558]{abram} of the hypergeometric function ${}_2 F_1(a,b;c;z)$. Now inserting the above result into formula $(\ref{func-F})$, we eventually deduce that
\begin{align*}
F(t) = & \frac{q\,\sin(\lambda\pi)\,\Gamma((\alpha+d-1)/2)\, \Gamma((d+2\lambda-3)/2)}{\pi^{(d+1)/2}\, \Gamma((d+\alpha+2\lambda-3)/2)}\, R_\ast^\alpha\, (R_\ast^2-t^2)^{\lambda-1}\, \bigg\{  - {}_2 F_1\left(-\frac{\alpha}{2}, 1;\lambda+1; 1-\left(\frac{t}{R_\ast}\right)^2\right)\\
 & \qquad \qquad \qquad + \frac{\alpha}{2\lambda(\lambda+1)}\bigg(1-\left(\frac{t}{R_\ast}\right)^2 \bigg)\, {}_2 F_1\left(1-\frac{\alpha}{2}, 2;\lambda+2; 1-\left(\frac{t}{R_\ast}\right)^2\right) \bigg\}, \quad 0\leq t\leq R_\ast.
\end{align*}

\noindent It remains to find the constant $C_Q$, which is computed according to $(\ref{constant-C})$. However, to avoid the tedious calculations encountered while evaluating the integral appearing on the right hand side of $(\ref{constant-C})$, we recall that the constant $C_Q$ is related to the Robin constant $F_Q$ by
\begin{equation*}
C_Q = \frac{\sin(\lambda\pi)\,\Gamma((d-1)/2)}{\pi^{(d+1)/2}}\, F_Q.
\end{equation*}
Furthermore, from the proof of Proposition \ref{ms-func-min-prop} it follows that $F_Q=\cF_s(\D_{R_\ast})$, where $R_\ast$ is given by $(\ref{mon-ext-field-extrem-meas-disk-rad})$. We now easily derive
\begin{equation*}
C_Q =   \frac {\Gamma((d+2\lambda-1)/2)} {\pi^{(d-1)/2}\, \Gamma(\lambda)} \, \bigg\{ \frac{1}{R_\ast^{d+2\lambda-3}} +  \frac{q\sin(\lambda\pi)\,\Gamma((d+\alpha-1)/2)\, \Gamma(\lambda)}{\pi\, \Gamma((d+\alpha+2\lambda-1)/2)}\, R_\ast^\alpha \bigg\}.
\end{equation*}

\qed

\noindent{\bf Proof of Theorem \ref{theo3}.} We commence by deriving expression $(\ref{pt-charge-ext-field-ext-mes})$ for the density of the extremal measure $\mu_Q$. For that upon inserting $(\ref{pt-charge})$ into $(\ref{func-g})$, we arrive at the following integral
\begin{align*}
\int_0^r \frac{u^{d-2}\,du}{(u^2+h^2)^{s/2}\, (r^2-u^2)^{1-\lambda}} 
& = \frac {r^{d+2\lambda-3}}{2(r^2+h^2)^{s/2}}\, \int_0^1 z^{\lambda-1}\, (1-z)^{(d-3)/2}\,(1-(r^2/(r^2+h^2))z)^{-s/2}\, dz \\
& = \frac{\Gamma((d-1)/2)\Gamma(\lambda)}{2\Gamma((d+2\lambda-1)/2)}\, \frac{r^{d+2\lambda-3}}{(r^2+h^2)^{s/2}}\, {}_2 F_1\left(\frac{s}{2}, \lambda; \frac{d+2\lambda-1}{2}; \frac{r^2}{r^2+h^2}\right)\\
& = \frac{\Gamma((d-1)/2)\Gamma(\lambda)}{2\Gamma((d+2\lambda-1)/2)}\, \frac{r^{d+2\lambda-3}}{(r^2+h^2)^{s/2}}\, {}_2 F_1\left(\frac{s}{2}, \lambda; \frac{s}{2}+1; \frac{r^2}{r^2+h^2}\right),
\end{align*}
where we used the substitution $r^2-u^2=r^2 z$ and the integral representation of the hypergeometric function \cite[\#15.3.1, p. 558]{abram}. Further recalling the known relation between the hypergeometric function ${}_2 F_1(a,b;c;z)$ and the incomplete Beta function $\Beta(z;a,b)$,
\begin{equation}\label{hypergeom-to-beta-1}
\Beta(z;a,b) = \frac{z^a}{a}\, {}_2 F_1(a,1-b;a+1;z),
\end{equation}
we eventually find that
\begin{equation}\label{pt-charge-int-1}
\int_0^r \frac{u^{d-2}\,du}{(u^2+h^2)^{s/2}\, (r^2-u^2)^{1-\lambda}} = \frac{(d+2\lambda-3)\, \Gamma(\lambda)\, \Gamma((d-1)/2)}{4\Gamma((d+2\lambda-1)/2)}\, \Beta\left(\frac{r^2}{r^2+h^2}; \frac{d+2\lambda-3}{2},1-\lambda\right).
\end{equation}
Differentiating both sides of $(\ref{pt-charge-int-1})$ with respect to $r$, while keeping in mind the fact
\begin{equation}\label{derivative-beta-func}
\frac{d}{dz}\, \Beta(z;a,b) = (1-z)^{b-1}\, z^{a-1},
\end{equation}
we deduce a concise expression for the function $g$, 
\begin{equation}\label{pt-charge-func-g}
g(r) = \frac{q (d+2\lambda-3)\, \Gamma(\lambda)\, \Gamma((d-1)/2)}{2\,\Gamma((d+2\lambda-1)/2)}\, \frac{h^{2(1-\lambda)}}{(r^2+h^2)^{(d-1)/2}}, \quad 0 \leq r \leq 1.
\end{equation}
Having expression for the function $g$ obtained, from formula $(\ref{func-F})$ it follows that we need to evaluate the integral
\begin{align*}
\int_t^1 \frac {r\,dr} {(r^2+h^2)^{(d-1)/2}\, (r^2-t^2)^{1-\lambda}} & = \frac {(1-t^2)^\lambda}{2}\, \int_0^1 \frac {z^{\lambda-1}\, dz}{(t^2+h^2+(1-t^2)z)^{(d-1)/2}} \\
	& = \frac {(1-t^2)^\lambda}{2\lambda\, (t^2+h^2)^{(d-1)/2}}\,  {}_2 F_1\left(\lambda, \frac{d-1}{2};\lambda+1; - \frac{1-t^2}{h^2+t^2}\right) \\
	& =  \frac {(1-t^2)^\lambda}{2\lambda\, (1+h^2)^{(d-1)/2}}\,  {}_2 F_1\left(1, \frac{d-1}{2};\lambda+1; \frac{1-t^2}{1+h^2}\right),
\end{align*}
where in the last two steps we used the integral representation of the hypergeometric function \cite[\#15.3.1, p. 558]{abram} and the linear transformation formula  \cite[\#15.3.5, p. 559]{abram}. Using another relation between the hypergeometric function ${}_2 F_1(a,b;c;z)$ and the incomplete Beta function $\Beta(z;a,b)$,
\begin{equation}\label{hypergeom-to-beta-2}
\Beta(z;a,b) = \frac{z^a\, (1-z)^b}{a}\, {}_2 F_1(1,a+b;a+1;z),
\end{equation}
we eventually obtain, after some trivial simplifications,
\begin{equation}\label{pt-charge-int-2}
\int_t^1 \frac {r\,dr} {(r^2+h^2)^{(d-1)/2}\, (r^2-t^2)^{1-\lambda}} = \frac{1}{2\,(h^2+t^2)^{(d-2\lambda-1)/2}}\, \Beta\left( \frac{1-t^2}{1+h^2}; \lambda, \frac{d-2\lambda-1}{2}\right).
\end{equation}
Therefore,
\begin{equation}\label{pt-charge-int-3}
\int_t^1 \frac{g(r) r dr}{(r^2-t^2)^{1-\lambda}} =  \frac{q (d+2\lambda-3)\, \Gamma(\lambda)\, \Gamma((d-1)/2)}{4\,\Gamma((d+2\lambda-1)/2)}\, \frac{h^{2(1-\lambda)}}{(h^2+t^2)^{(d-2\lambda-1)/2}}\, \Beta\left( \frac{1-t^2}{1+h^2}; \lambda, \frac{d-2\lambda-1}{2}\right).
\end{equation}
Differentiating both sides of $(\ref{pt-charge-int-3})$ with respect to $t$, and substituting the result into $(\ref{func-F})$, after some simple algebra we find
\begin{empheq}{align}\label{pt-charge-func-F}
F(t) = - \frac{q\,\sin(\lambda\pi)\,\Gamma((d-1)/2)}{\pi^{(d+1)/2}}\, h^{2(1-\lambda)} \, & \bigg\{ \frac{d-2\lambda-1}{2}\,  \frac{1}{(h^2+t^2)^{(d-2\lambda+1)/2}}\, \Beta\left( \frac{1-t^2}{1+h^2}; \lambda, \frac{d-2\lambda-1}{2}\right) \\ \nonumber
& + \frac{(1-t^2)^{\lambda-1}} {(1+h^2)^{(d-3)/2}\, (h^2+t^2)} \bigg\}, \quad 0\leq t \leq 1.
\end{empheq}
Our next step is to evaluate the constant $C_Q$, following the recipe contained in formula $(\ref{constant-C})$. Inserting expression $(\ref{pt-charge-func-F})$ into $(\ref{constant-C})$, we find that one of the integrals we need to evaluate is
\begin{align*}
\int_0^1 \frac{(1-t^2)^{\lambda-1}}{h^2+t^2}\, t^{d-2}\, dt & = \frac{1}{2h^2}\, \int_0^1 u^{(d-3)/2}\, (1-u)^{\lambda-1}\, (1-(-h^{-2})u)^{-1}\,du \\
	& = \frac{1}{2h^2}\, \frac{\Gamma(\lambda)\, \Gamma((d-1)/2)}{\Gamma((d+2\lambda-1)/2)}\, {}_2 F_1\left(1, \frac{d-1}{2}; \frac{d+2\lambda-1}{2}; -\frac{1}{h^2}\right)\\
	& = \frac{\Gamma(\lambda)\, \Gamma((d-1)/2)}{2\,\Gamma((d+2\lambda-1)/2)}\, \frac{1}{1+h^2}\, {}_2 F_1\left(1, \lambda; \frac{d+2\lambda-1}{2}; \frac{1}{1+h^2}\right),
\end{align*}
where the last two expressions were obtained using the integral representation of the hypergeometric function \cite[\#15.3.1, p. 558]{abram} and the linear transformation formula  \cite[\#15.3.4, p. 559]{abram}. We subsequently obtain
\begin{equation}\label{pt-charge-int-4}
\int_0^1 \frac{(1-t^2)^{\lambda-1}}{h^2+t^2}\, t^{d-2}\, dt  =  \frac{\Gamma(\lambda)\, \Gamma((d-1)/2)}{2\,\Gamma((d+2\lambda-1)/2)}\, \frac{1}{1+h^2}\, {}_2 F_1\left(1, \lambda; \frac{d+2\lambda-1}{2}; \frac{1}{1+h^2}\right).
\end{equation}
In light of $(\ref{pt-charge-int-4})$, we finally derive that
\begin{equation*}
C_Q = \frac {2\, \Gamma((d-1)/2+\lambda)} {\Gamma(\lambda)\, \Gamma((d-1)/2)}\, \bigg\{ \frac {\Gamma((d-1)/2)} {2\, \pi^{(d-1)/2}} + q\, c_{d,\lambda} \bigg\},
\end{equation*}
where
\begin{align*}
c_{d,\lambda}:= \frac{\sin(\lambda\pi)\,\Gamma((d-1)/2)}{\pi^{(d+1)/2}}\, h^{2(1-\lambda)} \, & \bigg\{ \frac{d-2\lambda-1}{2}\, \int_0^1 \frac{t^{d-2}}{(h^2+t^2)^{(d-2\lambda+1)/2}}\, \Beta\left( \frac{1-t^2}{1+h^2}; \lambda, \frac{d-2\lambda-1}{2}\right)\, dt \\ \nonumber
& + \frac {\Gamma(\lambda)\, \Gamma((d-1)/2)}  {2\, \Gamma((d-1)/2+\lambda)}  \, \frac{1}{(1+h^2)^{(d-1)/2}}\, {}_2 F_1\left(1, \lambda; \frac{d+2\lambda-1}{2}; \frac{1}{1+h^2}\right) \bigg\}.
\end{align*}

\noindent Now we show that the obtained extremal measure is indeed a positive measure, provided the height of the point charge satisfies certain restrictions. The first step toward that goal is to demonstrate that the density $f(r)$ of the extremal measure $\mu_Q$ is an increasing function of $r$. Differentiating its expression $(\ref{pt-charge-ext-field-dens})$ and dropping a positive term involving the incomplete Beta function, we see that
\begin{align*}
f'(r) & \geq C_Q \, \frac {2r(1-\lambda)}{(1-r^2)^{2-\lambda}} + q\frac{\sin(\lambda\pi)\Gamma((d-1)/2)}{\pi^{(d+1)/2}}\, \frac{h^{2(1-\lambda)}\, r\, (d-2\lambda+1)}{(h^2+r^2)^2\,(1-r^2)^{1-\lambda}\, (1+h^2)^{(d-3)/2}} \\
& \quad\qquad\qquad\qquad\qquad - q\frac{\sin(\lambda\pi)\Gamma((d-1)/2)}{\pi^{(d+1)/2}}\, \frac{h^{2(1-\lambda)}\, 2r\,(1-\lambda)}{(h^2+r^2)\,(1-r^2)^{2-\lambda}\, (1+h^2)^{(d-3)/2}} \\
& \geq \frac{\Gamma((d+2\lambda-1)/2)}{\pi^{(d-1)/2}\, \Gamma(\lambda)}\, \frac {2r(1-\lambda)}{(1-r^2)^{2-\lambda}} + q\frac{\sin(\lambda\pi)\Gamma((d-1)/2)}{\pi^{(d+1)/2}}\, \frac{h^{2(1-\lambda)}\, r\, (d-2\lambda+1)}{(h^2+r^2)^2\,(1-r^2)^{1-\lambda}\, (1+h^2)^{(d-3)/2}} \\
& \quad\qquad\qquad\qquad\qquad\qquad\qquad\qquad - q\frac{\sin(\lambda\pi)\Gamma((d-1)/2)}{\pi^{(d+1)/2}}\, \frac{h^{2(1-\lambda)}\, 2r\,(1-\lambda)}{(h^2+r^2)\,(1-r^2)^{2-\lambda}\, (1+h^2)^{(d-3)/2}}.
\end{align*}
Therefore,
\begin{align*}
& f'(r) \,\pi^{(d+1)/2}\,  (1-\lambda)^{-1}\,   (h^2+r^2)^2\, r^{-1}\, (1-r^2)^{2-\lambda}\,  (1+h^2)^{(d-3)/2} \\
& \geq \frac{ 2\pi \, \Gamma((d+2\lambda-1)/2)}{\Gamma(\lambda)} \, (h^2+r^2)^2\, (1+h^2)^{(d-3)/2} + q (d-2\lambda+1) (1-\lambda) \sin(\lambda\pi)\, \Gamma((d-1)/2)\, h^{2(1-\lambda)} \, (1-r^2)\\
& \qquad\qquad\qquad\qquad\qquad\qquad \qquad \qquad \qquad \qquad \qquad - q \sin(\lambda\pi) \Gamma((d-1)/2)\, h^{2(1-\lambda)} \, (h^2+r^2) \\
& \geq \frac{2  h^{d-3}\,  \Gamma((d+2\lambda-1)/2)} {\Gamma(\lambda)} \, \bigg\{ \pi\, (h^2+r^2)^2 + q\, \frac{(d-2\lambda+1)\,(1-\lambda)\, \sin(\lambda\pi)} {2}\, \Beta\left(\lambda, \frac{d-1}{2}\right) \, h^{5-2\lambda-d}\, (1-r^2) \\
& \qquad\qquad\qquad\qquad \qquad \qquad \qquad \qquad \qquad - q\, \frac{\sin(\lambda\pi)}{2}\, \Beta\left(\lambda, \frac{d-1}{2}\right) \, h^{5-2\lambda-d}\, (h^2+r^2) \bigg\}.
\end{align*}
Let $t=r^2+h^2$. Then the expression in the braces in the last line above becomes
\begin{align*}
 \pi\, (h^2+r^2)^2 & + q\, \frac{(d-2\lambda+1)\,(1-\lambda)\, \sin(\lambda\pi)} {2}\, \Beta\left(\lambda, \frac{d-1}{2}\right) \, h^{5-2\lambda-d}\, (1-r^2) \\
& - q\, \frac{(1-\lambda)\,\sin(\lambda\pi)}{2}\, \Beta\left(\lambda, \frac{d-1}{2}\right) \, h^{5-2\lambda-d}\, (h^2+r^2) \\
& = \pi t^2 - q\, \frac{((1-\lambda)(d-2\lambda+1)+1) \sin(\lambda\pi)}{2} \Beta\left(\lambda, \frac{d-1}{2}\right)\, h^{5-2\lambda-d}\, t \\
&   + q\, \frac{(d-2\lambda+1)\,(1-\lambda)\,\sin(\lambda\pi)}{2}\, \Beta\left(\lambda, \frac{d-1}{2}\right)\, h^{5-2\lambda-d}\, (1+h^2).
\end{align*}
The quadratic polynomial 
\begin{align*}
m(t):= \pi t^2 & - q\, \frac{((1-\lambda)(d-2\lambda+1)+1) \sin(\lambda\pi)}{2} \Beta\left(\lambda, \frac{d-1}{2}\right)\, h^{5-2\lambda-d}\, t \\
	& + q\, \frac{(d-2\lambda+1)\,(1-\lambda)\,\sin(\lambda\pi)}{2}\, \Beta\left(\lambda, \frac{d-1}{2}\right)\, h^{5-2\lambda-d}\, (1+h^2)
\end{align*}
has the discriminant
\begin{align*}
D = - \frac{q\sin(\lambda\pi) h^{2(5-2\lambda-d)}}{4} \,  \Beta  \left(  \lambda, \frac{d-1}{2}\right)\, \bigg\{ & 8\pi(d+2\lambda-1)(1-\lambda)\, (1+h^2)\, h^{d+2\lambda-5} \\
- & q ((1-\lambda)(d-2\lambda+1)+1)^2 \sin(\lambda\pi)\, \Beta\left(\lambda, \frac{d-1}{2}\right) \bigg\}.
\end{align*}
For the term in braces in the last expression we have the following trivial estimate
\begin{align*}
& 8\pi(d+2\lambda-1)(1-\lambda)\, (1+h^2)\, h^{d+2\lambda-5} - q ((1-\lambda)(d-2\lambda+1)+1)^2 \sin(\lambda\pi)\, \Beta\left(\lambda, \frac{d-1}{2}\right) \\
& \geq 8\pi(d+2\lambda-1)(1-\lambda)\,  h^{d+2\lambda-3}  - q ((1-\lambda)(d-2\lambda+1)+1)^2 \sin(\lambda\pi)\, \Beta\left(\lambda, \frac{d-1}{2}\right) \\
& > 0,
\end{align*}
if $h$ is chosen so that $h>h_-$,
where
\begin{equation*}
h_-:=\left(\frac{q\, ((1-\lambda)(d-2\lambda+1)+1)^2 \sin(\lambda\pi)}{8\pi(d+2\lambda-1)(1-\lambda)} \, \Beta\left(\lambda, \frac{d-1}{2}\right) \right)^{1/(d+2\lambda-3)}
\end{equation*}
With such a choice of $h$ it is clear that the discriminant $D$ is strictly negative, and hence the polynomial $m(t)$ does not have real roots. Therefore, $m(t)>0$, which in turn shows that $f'(r)>0$. Thus, $f(r)>f(0)$, for all $0<r<1$ and if $h$ chosen such that $h>h_-$.

Set $p(h):=f(0)$. We then obtain
\begin{align*}
p(h) & = \frac {2\, \Gamma((d-1)/2+\lambda)} {\Gamma(\lambda)\, \Gamma((d-1)/2)}\, \bigg\{ \frac {\Gamma((d-1)/2)} {2\, \pi^{(d-1)/2}} + q\, c_{d,\lambda} \bigg\}  \\
& - \frac {q\,\sin(\lambda\pi)\,\Gamma((d-1)/2)} {\pi^{(d+1)/2}} \,  \bigg\{ \frac{d-2\lambda-1}{2}\,  \frac{1} {h^{d-1}}\, \Beta\left( \frac{1}{1+h^2}; \lambda, \frac{d-2\lambda-1}{2}\right)  + \frac{1} {h^{2\lambda} (1+h^2)^{(d-3)/2}} \bigg\}.
\end{align*}
It is not hard to see that
\begin{equation*}
\lim_{h\to0+} p(h) = -\infty, \quad \lim_{h\to\infty} p(h) = \frac {\Gamma((d-1)/2+\lambda)} {\pi^{(d-1)/2}\, \Gamma(\lambda)} > 0.
\end{equation*}
As $p(h)$ is a continuous function, it has at least one positive root. Denote the largest such root by $h_+$. It then follows that $f(r)>0$ provided $h>\max\{h_-,h_+\}$, and therefore $\mu_Q$ is a positive measure, as required.

\qed


\noindent{\bf Proof of Corollary \ref{pnt-charge-newtonian-case}.}  We begin by evaluating the integral appearing in the right hand side of the constant $c_{d,q}$,
\begin{empheq}{align}\label{complicated-int-1}
 \int_0^1 \frac{t^{d-2}}{(h^2+t^2)^{(d-2\lambda+1)/2}}\, & \Beta\left( \frac{1-t^2}{1+h^2}; \lambda, \frac{d-2\lambda-1}{2}\right)\, dt \\ \nonumber
 & = \frac{1}{2\lambda}\, \xi^{(d+1)/2}\, \int_0^1 u^\lambda\, (1-u)^{(d-3)/2}\, (1-u\xi)^{-1}\,  {}_2 F_1\left(1,\frac{d-1}{2};\lambda+1;u\xi\right),
\end{empheq}
where we put $\xi=1/(1+h^2)$ for brevity. Recall that in the case of Newtonian potential we have $s=d-2$, where $d=2m+4$, with $m\geq2$ is a natural number. Using a well-known fact  \cite[\#15.4.1, p. 561]{abram} that the hypergeometric function $ {}_2 F_1(a,b;c;z)$ reduces to a polynomial if either one of its first two parameters is a negative integer, we easily find
\begin{equation}\label{hypergeom-to-poly}
{}_2 F_1\left(1,\frac{d-1}{2};\frac{3}{2};u\xi\right) = (1-u\xi)^{-(m+1)}\, \sum_{n=0}^m \frac {(-m)_n}{(2n+1)n!}\, u^n \xi^n.
\end{equation}
Inserting expression $(\ref{hypergeom-to-poly})$ into the integral on the right hand side of $(\ref{complicated-int-1})$ and again using  \cite[\#15.4.1, p. 561]{abram}, we eventually find that
\begin{empheq}{align}\label{complic-int-to-elem-func}
& \int_0^1 \frac{t^{d-2}}{(h^2+t^2)^{(d-2\lambda+1)/2}}\,  \Beta\left( \frac{1-t^2}{1+h^2}; \lambda, \frac{d-2\lambda-1}{2}\right)\, dt \\ \nonumber
 &  = \Gamma\left(m+\frac{3}{2}\right)\, \xi^{m+5/2}\, \sum_{n=0}^m \frac{(-m)_n}{(2n+1)n!}\,\xi^n \, \sum_{l=0}^{m-2} \frac{(2-m)_l\,\Gamma(n+l+3/2)}{(n+m+l+2)!}\, \xi^l \\\nonumber
 & = \Gamma\left(m+\frac{3}{2}\right)\, (1+h^2)^{-(m+5/2)}\, \sum_{n=0}^m \frac{(-m)_n}{(2n+1)n!}\,(1+h^2)^{-n} \, \sum_{l=0}^{m-2} \frac{(2-m)_l\,\Gamma(n+l+3/2)}{(n+m+l+2)!}\,(1+h^2)^{-l}.
 \end{empheq}
 We also note that the hypergeometric function appearing in the definition of the constant $c_{d,\lambda}$ in this special case reduces to a finite sum involving elementary functions only. Indeed, using  \cite[\#15.3.19, p. 560]{abram}, with some work we can show that
\begin{equation}\label{const-c-d-lambda-heypergeom-to-elem}
{}_2 F_1\left(1, \lambda; \frac{d+2\lambda-1}{2}; \frac{1}{1+h^2}\right) = 2(m+1)! \frac{\sqrt{1+h^2}}{h+\sqrt{1+h^2}}\, \sum_{n=0}^m \frac{(-m)_n}{(m+n+1)!}\, \left(\frac{\sqrt{1+h^2}-h}{\sqrt{1+h^2}+h}\right)^n.
\end{equation}
We also reduce the Beta-function term, appearing in expression $(\ref{pt-charge-func-F-theo})$, to a finite sum of expressions in elementary functions only. For that we first recall the expression
\begin{equation*}
\Beta\left( \frac{1-r^2}{1+h^2}; \lambda, \frac{d-2\lambda-1}{2}\right) = \frac{1}{\lambda}\, \frac{(1-r^2)^\lambda}{(1+h^2)^{(d-1)/2}}\, (h^2+r^2)^{(d-2\lambda-1)/2}\, {}_2 F_1\left(1,\frac{d-1}{2};\lambda+1;\frac{1-r^2}{1+h^2} \right),
\end{equation*}
which, after taking into account that $\lambda=1/2$ and $d=2m+4$, $m\geq 2$, and using the fact  \cite[\#15.4.1, p. 561]{abram}, gives us
\begin{empheq}{align}\label{func-F-beta-to-elem}
\frac{d-2\lambda-1}{2}\, \frac{1}{(h^2+r^2)^{(d-2\lambda+1)/2}} & \, \Beta\left( \frac{1-r^2}{1+h^2}; \lambda, \frac{d-2\lambda-1}{2}\right) \\ \nonumber
& =\frac{ 2(m+1)}{(h^2+r^2)^{m+2}}\,  \sum_{n=0}^m \frac{(-m)_n}{(2n+1)n!}\, \left( \frac{1-r^2}{1+h^2}\right)^{n+1/2}.
\end{empheq}
With the help of $(\ref{func-F-beta-to-elem})$, expression $(\ref{pt-charge-func-F-theo})$ is reduced to
\begin{empheq}{align*}
F(r)= - qh \frac{\Gamma(m+3/2)}{\pi^{m+5/2}}\,\bigg\{ & \frac{ 2(m+1)}{(h^2+r^2)^{m+2}}\,   \sum_{n=0}^m \frac{(-m)_n}{(2n+1)n!}\, \left( \frac{1-r^2}{1+h^2}\right)^{n+1/2} \\ \nonumber
& + \frac{1}{\sqrt{1-r^2}}\, \frac{1}{(h^2+r^2)\,(1+h^2)^{m+1/2}} \bigg\}, \quad 0\leq r\leq 1.
\end{empheq}
We similarly derive that in the case $\lambda=1/2$ and $d=2m+4$, $m\geq 2$,
\begin{align*}
 c_{d,\lambda} & = \frac{(\Gamma(m+3/2))^2}{\pi^{m+5/2}}\, h\, (1+h^2)^{-(m+1)} \\
&  \bigg\{  \Gamma\left(m+\frac{3}{2}\right)\, (1+h^2)^{-(m+5/2)}\, \sum_{n=0}^m \frac{(-m)_n}{(2n+1)n!}\,(1+h^2)^{-n} \, \sum_{l=0}^{m-2} \frac{(2-m)_l\,\Gamma(n+l+3/2)}{(n+m+l+2)!}\,(1+h^2)^{-l} \\ \nonumber
& + \frac{\sqrt{1+h^2}}{h+\sqrt{1+h^2}}\, \sum_{n=0}^m \frac{(-m)_n}{(m+n+1)!}\, \left(\frac{\sqrt{1+h^2}-h}{\sqrt{1+h^2}+h}\right)^n \bigg\},
\end{align*}
and also 
\begin{equation*}
C_Q = \frac{2(m+1)!}{\sqrt{\pi}\Gamma(m+3/2)}\,\bigg\{ \frac{\Gamma(m+3/2)}{2\pi^{m+3/2}} + q\, c_{d,\lambda} \bigg\}.
\end{equation*}
\qed


\noindent{\bf Proof of Corollary \ref{disk-three-dim-coulomb-pot}.} The expression $(\ref{density-disk-three-dim-coulomb-pot})$ for the density $f(r)$ of the extremal measure $\mu_Q$ can be derived either using the recipe given in Theorem \ref{theo1}, or alternatively from corresponding parts of Theorem \ref{mon-ext-field-extrem-meas}.

After obtaining the required expression for the density $f(r)$, the next step is to show sure that the obtained solution is, in fact, represents a positive measure. Below we will show that this is indeed the case, provided a point charge is located sufficiently far from the surface of the disk $\D$. In the course of the proof we precisely determine the critical hight of location of the point charge, guaranteeing the positivity of the extremal measure. 

We begin by observing that the density $f(r)$ is a strictly increasing function of $r$. Differentiating expression $(\ref{density-disk-three-dim-coulomb-pot})$, after simplifications, one finds that
\begin{empheq}{align*}
 f'(r) = & \frac{c}{\pi^2} \frac{r}{(1-r^2)^{3/2}}+\frac{3h}{\pi^2}\frac{1}{\sqrt{1-r^2}} \frac{r}{(h^2+r^2)^2}\\
 	& - \frac{h}{\pi^2} \frac{1}{(1-r^2)^{3/2}} \frac{r}{h^2+r^2} + \frac{3h}{\pi^2}\frac{r}{(h^2+r^2)^{5/2}}\tan^{-1}\sqrt{\frac{1-r^2}{h^2+r^2}},
\end{empheq}
with
\begin{equation*}
c=\frac{\pi}{2}\left(1+\frac{2h\tan^{-1}(1/h)}{\pi\sqrt{1+h^2}}\right).
\end{equation*}
Using trivial estimates, we further obtain
\begin{empheq}{align*}
\pi^2f'(r)(1-r^2)^{3/2}(h^2+r^2)^2r^{-1} & \geq c (h^2+r^2)^2+3h(1-r^2) - h(h^2+r^2)\\
					& \geq \frac{\pi}{2} (h^2+r^2)^2 +3h(1-r^2) - h(h^2+r^2)\\
					& = \frac{1}{2} (\pi t^2-8ht+6h(1+h^2)),
\end{empheq}
with $t=r^2+h^2$. The quadratic polynomial $m(t)=\pi t^2-8ht+6h(1+h^2)$ has discriminant $D=-8h(3\pi(1+h^2)-8h)<0$. Hence $m(t)$ does not have real roots, which implies that $m(t)>0$. This in turn shows that $f'(r)>0$, as required, and so $f(r)>f(0)$ for all $0<r<1$. 

Denoting $p(h):=f(0)$, we readily find
\begin{equation*}
p(h) = \frac{1}{2\pi}\left(1+\frac{2h\tan^{-1} (1/h)}{\pi\sqrt{1+h^2}} \right) - \frac{1}{\pi^2 h} - \frac{1}{\pi^2h^2}\tan^{-1} (1/h).
\end{equation*}
As observed during the course of the proof of Theorem \ref{mon-ext-field-extrem-meas}, we have
\begin{equation*}
\lim_{h\to0+} p(h) = -\infty, \quad \lim_{h\to\infty} p(h) = \frac{1}{2\pi}  > 0,
\end{equation*}
and so as $p(h)$ is a continuous function, it has at least one positive root. Denote the largest such a root by $h_+$. We will show that $h_+$ is the unique root of $p(h)$ on $[0,\infty)$. This follows from the fact that $p(h)$ is a strictly increasing function on $[0,\infty)$, which we now demonstrate. Indeed, the derivative of $p(h)$ can be easily computed to be
\begin{empheq}{align}\label{func-p-derivative}
p'(h) & = \frac{\tan^{-1}(1/h)}{\pi^2\sqrt{1+h^2}} -  \frac{h}{\pi^2(1+h^2)^{3/2}} - \frac{h^2\tan^{-1}(1/h)}{\pi^2(1+h^2)^{3/2}} \\\nonumber
	& + \frac{1}{\pi^2 h^2} + \frac{2\tan^{-1}(1/h)}{\pi^2 h^3} + \frac{1}{\pi^2 h^2 (1+h^2)}.
\end{empheq}
From the Mean Value Theorem for the derivative it follows that for $x>0$,
\begin{equation}\label{tan-inequal}
\frac{x}{1+x^2} < \tan^{-1}x < x.
\end{equation}
Using inequality $(\ref{tan-inequal})$ in conjunction with $(\ref{func-p-derivative})$, we write down an estimate of $p'(h)$ from below,
\begin{align*}
p'(h) & \geq \frac{h}{\pi^2(1+h^2)^{3/2}} - \frac{h}{\pi^2(1+h^2)^{3/2}} - \frac{h}{\pi^2(1+h^2)^{3/2}} \\
	& + \frac{1}{\pi^2 h^2} +  \frac{2}{\pi^2 h^2(1+h^2)} + \frac{1}{\pi^2 h^2(1+h^2)} \\
	& =  \frac{3}{\pi^2 h(1+h^2)} + \frac{1}{\pi^2 h^2} - \frac{h}{\pi^2(1+h^2)^{3/2}} \\
	& \geq \frac{3}{\pi^2 h^2(1+h^2)} + \frac{1}{\pi^2 h^2} - \frac{1}{\pi^2 h^2} \\
	& = \frac{3}{\pi^2 h^2(1+h^2)} > 0, \quad h > 0.
\end{align*}
Clearly, $p(h)\geq 0$ when $h\geq h_+$ and $p(h)<0$ when $h<h_+$. This implies that $h=h_+$ is the critical height of the point charge, and so $f(r)\geq0$ for all $0\leq r\leq 1$, provided $h\geq h_+$, which means that $\mu_Q$ is a positive measure.

If $h<h_+$ then support of the extremal measure will no longer be the entire disk, but rather will have an opening around the origin. Indeed, if this is not so, from the fact that $p(h)<0$ when $h<h_+$, we see that in this case the density $f(r)$ will be negative for all $0\leq r \leq 1$. This implies that the measure $\mu_Q$ in $(\ref{mes-disk-three-dim-coulomb-pot})$ will no longer be a positive measure, and thus cannot be the extremal measure. This contradiction shows that when $h<h_+$ the point charge will clear out an opening in the disk $\D$ at the origin.

\qed


\noindent{\bf Proof of Theorem \ref{theo5}.} Recalling that the extremal measure is absolutely continuous with respect to the Lebesgue surface area measure, and invoking Lemma $\ref{integ-rep-form}$, for the Riesz $s$-potential $U_s^{\mu_Q}$ we have the following representation, for $x=(0,r\overbar{x}) \in\D$, with $\overbar{x}\in\S^{d-2}$ and $0\leq a\leq r \leq b\leq 1$,
\begin{empheq}{align*}\label{newtpot-1}
U_s^{\mu_Q}(x) & = \int\frac{1}{|x-y|^s}\,d\mu(y)\\ \nonumber
			& = \frac {2\pi^{(d-2)/2}} {\Gamma(d/2-1)}\, \int_a^b f(\rho)\, \rho^{d-2}\, d\rho \, \int_0^\pi \frac {\sin^{d-3}\xi \,d\xi} {(r^2 + \rho^2 - 2r\rho \cos\xi)^{s/2}}\\ \nonumber
			& = \frac {4\sin(\lambda\pi)\,\pi^{(d-3)/2}\,\Gamma(\lambda)} {\Gamma((d+2\lambda-3)/2)}\, r^{3-d}\,  \int_a^b f(\rho)\, \rho \, d\rho \, \int_0^{\min{(r,\rho)}} \frac {t^{d+2\lambda-4}\,dt} {(r^2-t^2)^\lambda \, (\rho^2-t^2)^\lambda} 
\end{empheq}
This allows us to write integral equation $(\ref{ieec})$ in the following form,
\begin{equation}\label{ie1-ch-5}
\int_a^b \int_0^{\min(\rho,r)} \frac{f(\rho)\,\rho\, t^{d+2\lambda-4}\, dt\,d\rho} {(r^2-t^2)^\lambda \, (\rho^2-t^2)^\lambda} = \frac  {\Gamma((d+2\lambda-3)/2)} {4\sin(\lambda\pi)\,\pi^{(d-3)/2}\,\Gamma(\lambda)}\, r^{d-3} \left(F_Q - Q(r)\right), \quad a\leq r \leq b.
\end{equation}
We continue by working with the left hand side of equation $(\ref{ie1-ch-5})$. Splitting the domain of integration for the variable $\rho$, we write
\begin{empheq}{align}\label{int-0-ch-5}
\int_a^b \int_0^{\min(\rho,r)} \frac{f(\rho)\,\rho\, t^{d+2\lambda-4}\, dt\,d\rho} {(r^2-t^2)^\lambda \, (\rho^2-t^2)^\lambda}
& = \int_a^r \int_0^\rho\frac{f(\rho)\,\rho\, t^{d+2\lambda-4}\, dt\,d\rho} {(r^2-t^2)^\lambda \, (\rho^2-t^2)^\lambda} \\\nonumber
& + \int_r^b\int_0^r \frac{f(\rho)\,\rho\, t^{d+2\lambda-4}\, dt\,d\rho} {(r^2-t^2)^\lambda \, (\rho^2-t^2)^\lambda}.
\end{empheq}
In the first integral on the right hand side of the above expression, we further split the domain of integration as follows,
\begin{empheq}{align}\label{int-1-ch-5}
\int_a^r \int_0^\rho\frac{f(\rho)\,\rho\, t^{d+2\lambda-4}\, dt\,d\rho} {(r^2-t^2)^\lambda \, (\rho^2-t^2)^\lambda} & 
		= \int_a^r \int_0^a \frac{f(\rho)\,\rho\, t^{d+2\lambda-4}\, dt\,d\rho} {(r^2-t^2)^\lambda \, (\rho^2-t^2)^\lambda} \\\nonumber
		& + \int_a^r \int_a^\rho \frac{f(\rho)\,\rho\, t^{d+2\lambda-4}\, dt\,d\rho} {(r^2-t^2)^\lambda \, (\rho^2-t^2)^\lambda}.
\end{empheq}
Similarly, the second integral is split in the following way,
\begin{empheq}{align}\label{int-2-ch-5}
\int_r^b \int_0^r \frac{f(\rho)\,\rho\, t^{d+2\lambda-4}\, dt\,d\rho} {(r^2-t^2)^\lambda \, (\rho^2-t^2)^\lambda} & 
		= \int_r^b\int_0^a \frac{f(\rho)\,\rho\, t^{d+2\lambda-4}\, dt\,d\rho} {(r^2-t^2)^\lambda \, (\rho^2-t^2)^\lambda} \\ \nonumber
		& + \int_r^b\int_a^r \frac{f(\rho)\,\rho\, t^{d+2\lambda-4}\, dt\,d\rho} {(r^2-t^2)^\lambda \, (\rho^2-t^2)^\lambda}.
\end{empheq}
We then change the order of integration in the second integral on the right hand side of $(\ref{int-1-ch-5})$ as follows,
\begin{equation}\label{int-3-ch-5}
\int_a^r \int_a^\rho\frac{f(\rho)\,\rho\, t^{d+2\lambda-4}\, dt\,d\rho} {(r^2-t^2)^\lambda \, (\rho^2-t^2)^\lambda} = \int_a^r \int_t^r \frac{f(\rho)\,\rho\, t^{d+2\lambda-4}\, dt\,d\rho} {(r^2-t^2)^\lambda \, (\rho^2-t^2)^\lambda}.
\end{equation}
Combining the first integral on the right hand side of $(\ref{int-1-ch-5})$ with the first integral on the right hand side of $(\ref{int-2-ch-5})$, we obtain 
\begin{empheq}{align}\label{int-4-ch-5}
\int_a^r \int_0^a \frac{f(\rho)\,\rho\, t^{d+2\lambda-4}\, dt\,d\rho} {(r^2-t^2)^\lambda \, (\rho^2-t^2)^\lambda} & + \int_r^b \int_0^a \frac{f(\rho)\,\rho\, t^{d+2\lambda-4}\, dt\,d\rho} {(r^2-t^2)^\lambda \, (\rho^2-t^2)^\lambda} \\ \nonumber
& = \int_a^b f(\rho)\,\rho\,d\rho\, \int_0^a \frac{t^{d+2\lambda-4}\, dt} {(r^2-t^2)^\lambda \, (\rho^2-t^2)^\lambda},
\end{empheq}
while similarly combining the integral on the right hand side of $(\ref{int-3-ch-5})$ with the second integral on the right hand side of $(\ref{int-2-ch-5})$ yields
\begin{empheq}{align}\label{int-5-ch-5}
\int_a^r \int_t^\rho \frac{f(\rho)\,\rho\, t^{d+2\lambda-4}\, dt\,d\rho} {(r^2-t^2)^\lambda \, (\rho^2-t^2)^\lambda} & + \int_r^b \int_a^r \frac{f(\rho)\,\rho\, t^{d+2\lambda-4}\, dt\,d\rho} {(r^2-t^2)^\lambda \, (\rho^2-t^2)^\lambda} \\ \nonumber
& = \int_a^r \frac{t^{d+2\lambda-4}\, dt}{(r^2-t^2)^\lambda} \int_t^b \frac{f(\rho)\,\rho\,d\rho}{(\rho^2-t^2)^\lambda}.
\end{empheq}
Collecting the above calculations, we conclude that integral equation $(\ref{ie1-ch-5})$ is transformed into
\begin{empheq}{align}\label{int-6-ch-5}
 \int_a^r \frac{t^{d+2\lambda-4}\, dt}{(r^2-t^2)^\lambda}  \int_t^b \frac{f(\rho)\, \rho\,d\rho}{(\rho^2-t^2)^\lambda} & +  \int_a^b f(\rho)\,\rho\,d\rho\, \int_0^a \frac{t^{d+2\lambda-4}\, dt} {(r^2-t^2)^\lambda \, (\rho^2-t^2)^\lambda} \\ \nonumber
 &  =  \frac  {\Gamma((d+2\lambda-3)/2)} {4\sin(\lambda\pi)\,\pi^{(d-3)/2}\,\Gamma(\lambda)}\, r^{d-3} \left(F_Q - Q(r)\right),  \quad a\leq r \leq b.
\end{empheq}
Our goal now is to further transform the second term on the left hand side of equation $(\ref{int-6-ch-5})$. For that we introduce the function $g(u,t)$ as a solution of the following Abel-type integral equation,
\begin{equation}\label{func-g-ch-5}
\int_a^\rho \frac{g(u,t)\,du}{(\rho^2-u^2)^\lambda} = \frac{1}{(\rho^2-t^2)^\lambda}, \quad a \leq \rho \leq b,
\end{equation}
where the variable $t$, such that $0\leq t \leq a$, is fixed. We thus obtain
\begin{empheq}{align}\label{int-7-ch-5}
\int_a^b & f(\rho)\,\rho\, d\rho  \int_0^a  \frac{t^{d+2\lambda-4}\, dt} {(r^2-t^2)^\lambda \, (\rho^2-t^2)^\lambda} \\ \nonumber
& = \int_a^b f(\rho)\,\rho\, d\rho \int_0^a \frac{t^{d+2\lambda-4}\,dt}{(r^2-t^2)^\lambda} \int_a^\rho \frac{g(u,t)\,du}{ (\rho^2-t^2)^\lambda} \\\nonumber
&  = \int_a^b f(\rho)\,\rho\, d\rho \int_0^a t^{d+2\lambda-4}\, dt \int_a^r \frac{g(s,t)\,ds}{(r^2-s^2)^\lambda}  \int_a^\rho \frac{g(u,t)\,du}{(\rho^2-u^2)^\lambda} \\\nonumber
& =  \int_a^r \frac{ds}{(r^2-s^2)^\lambda} \bigg\{ \int_a^b f(\rho)\,\rho\, d\rho \int_a^\rho \frac{du}{(\rho^2-u^2)^\lambda}  \int_0^a g(s,t)\, g(u,t)\, t^{d+2\lambda-4}\, dt \bigg\} \\\nonumber
& =  \int_a^r \frac{ds}{(r^2-s^2)^\lambda} \bigg\{  \int_a^b du \int_u^b \frac{f(\rho)\,\rho\, d\rho}{(\rho^2-u^2)^\lambda} \, \bigg\{ \int_0^a g(s,t)\, g(u,t)\, t^{d+2\lambda-4}\, dt \bigg\} \bigg\}.
\end{empheq}
Let
\begin{equation}\label{func-G-ch-5}
G(s) = \int_s^b \frac {f(u)\, u\, du} {(u^2-s^2)^\lambda}.
\end{equation}
Combining $(\ref{int-7-ch-5})$ and $(\ref{func-G-ch-5})$, we recast equation $(\ref{int-6-ch-5})$ into
\begin{align*}
 \int_a^r \frac{G(t)\, dt}{(r^2-t^2)^\lambda}  + \int_a^r \frac{ds}{(r^2-s^2)^\lambda} & \bigg\{  \int_a^b G(u)\,du \bigg\{ \int_0^a g(s,t)\, g(u,t)\,  t^{d+2\lambda-4}\, dt \bigg\} \bigg\}  \\
 & = \frac  {\Gamma((d+2\lambda-3)/2)} {4\sin(\lambda\pi)\,\pi^{(d-3)/2}\,\Gamma(\lambda)}\, r^{d-3} \left(F_Q - Q(r)\right),  \quad a\leq r \leq b.
\end{align*}
or,
\begin{empheq}{align}\label{int-8-ch-5}
 \int_a^r \frac{ds}{(r^2-s^2)^\lambda} \bigg\{G(s) + \int_a^b G(u)\, du\, & \bigg\{ \int_0^a g(s,t)\, g(u,t)\, t^{d+2\lambda-4}\, dt \bigg\} \bigg\} \\ \nonumber
 		& =\frac  {\Gamma((d+2\lambda-3)/2)} {4\sin(\lambda\pi)\,\pi^{(d-3)/2}\,\Gamma(\lambda)}\, r^{d-3} \left(F_Q - Q(r)\right),  \quad a\leq r \leq b.
\end{empheq}
Equation $(\ref{int-8-ch-5})$ is an Abel-type integral equation with respect to the function
\[
G(s) + \int_a^b G(u)\,du\, \bigg\{ \int_0^a g(s,t)\, g(u,t)\,  t^{d+2\lambda-4}\, dt \bigg\}.
\]
Solving this equation  \cite[\# 44, p. 122]{polman}, we obtain the following integral equation,
\begin{empheq}{align}\label{int-9-ch-5}
G(r) + \int_a^b G(u)\, du\, & \bigg\{  \int_0^a g(r,t)\, g(u,t)\,  t^{d+2\lambda-4}\, dt  \bigg\} \\\nonumber 
& = \frac  {\Gamma((d+2\lambda-3)/2)} {2\pi^{(d-1)/2}\,\Gamma(\lambda)}\, \frac{d}{dr} \int_a^r \frac{(F_Q-Q(\rho))\,\, \rho^{d-3}\, \rho\,d\rho}{(r^2-\rho^2)^{1-\lambda}}, \quad a \leq r \leq b.
\end{empheq}
We now turn our attention to evaluating the following expression, present on the right hand side of last expression,
\begin{equation*}
 \frac{d}{dr} \int_a^r \frac {\rho^{d-3}\, \rho\,d\rho}{(r^2-\rho^2)^{1-\lambda}}.
\end{equation*}
Using the substitution $r^2-t^2=r^2 z$, after some elementary calculations we find that
\[
 \int_a^r \frac {\rho^{d-3}\, \rho\,d\rho}{(r^2-\rho^2)^{1-\lambda}} = \frac{1}{2} \, r^{d+2\lambda-3}\, \Beta \left(1-\left(\frac{a}{r}\right)^2;\lambda, \frac{d-1}{2}\right),
\]
where $\Beta(z;a,b)$ is the incomplete Beta function defined in $(\ref{betafdef})$.

Now it is easy to see that
\begin{equation}\label{ring-int-diff-1}
\frac{d}{dr} \int_a^r \frac {\rho^{d-3}\, \rho\,d\rho}{(r^2-\rho^2)^{1-\lambda}} =  \frac{d+2\lambda-3}{2} \, r^{d+2\lambda-4}\, \Beta \left(1-\left(\frac{a}{r}\right)^2;\lambda, \frac{d-1}{2}\right) + \frac {a^{d-1}}{r}\, (r^2-a^2)^{\lambda-1}.
\end{equation}

\noindent Finally we deal with the inner integral on the left hand side of equation $(\ref{int-9-ch-5})$. First we recover the function $g$ from integral equation $(\ref{func-g-ch-5})$. Applying \cite[\#41, p. 11]{polman}, we find that the function $g$ is given by
\begin{equation}\label{ring-function-g-eval}
g(\rho,t) = \frac{2\sin(\lambda\pi)}{\pi} \, \frac{d}{d\rho} \int_a^\rho \frac {u\,du}{(\rho^2-u^2)^{1-\lambda}\, (u^2-t^2)^\lambda}, \quad 0\leq t\leq  \rho \leq b.
\end{equation}
Using the substitution $\rho^2-u^2=z$, after some work we find that
\begin{equation*}
 \int_a^\rho \frac {u\, du}{(\rho^2-u^2)^{1-\lambda}\, (u^2-t^2)^\lambda} = \frac{1}{2\lambda} \bigg( \frac{\rho^2-a^2}{\rho^2-t^2} \bigg)^\lambda \, {}_2 F_1\left(\lambda,\lambda;\lambda+1; \frac{\rho^2-a^2}{\rho^2-t^2}\right).
\end{equation*}
Taking derivative of the last expression with respect to $\rho$, and inserting the result into $(\ref{ring-function-g-eval})$ produces a remarkably simple expression for the function $g$,
\begin{equation}\label{func-g-sol-ch-5}
g(\rho,t)=\frac{2\sin(\lambda\pi)}{\pi} \frac{\rho}{\rho^2-t^2} \, \bigg( \frac{\rho^2-a^2}{a^2-t^2} \bigg)^{\lambda-1}, \quad  0\leq t \leq a \leq \rho \leq b.
\end{equation}
We therefore obtain
\begin{empheq}{align*}
 \int_0^a  & g(r,t)  g(u,t)\,  t^{d+2\lambda-4}\, dt \\
  & = \left(\frac{2\sin(\lambda\pi)}{\pi}\right)^2  ru \, (r^2-a^2)^{\lambda-1}  (u^2-a^2)^{\lambda-1}
 \int_0^a \frac{(a^2-t^2)^{2(1-\lambda)} \, t^{d+2\lambda-4}\, dt } {(r^2-t^2)\, (u^2-t^2)}.
\end{empheq}
After simple but fairly laborious calculations one finds that
\begin{align*}
& \int_0^a \frac{ (a^2-t^2)^{2(1-\lambda)} \, t^{d+2\lambda-4}\, dt } {(r^2-t^2)\, (u^2-t^2)} = \frac {\Gamma((d+2\lambda-3)/2)\Gamma(3-2\lambda)}{2\,\Gamma((d-2\lambda+3)/2)} \, \frac{a^{d-2\lambda+1}}{u^2-r^2} \\
& \times \bigg\{ \frac{1}{r^2}\,  {}_2 F_1\left(1,\frac{d+2\lambda-3}{2};\frac{d-2\lambda+3}{2};\left(\frac{a}{r}\right)^2\right) 
- \frac{1}{u^2}\,  {}_2 F_1\left(1,\frac{d+2\lambda-3}{2};\frac{d-2\lambda+3}{2};\left(\frac{a}{u}\right)^2\right) \bigg\}.
\end{align*}
Denote
\begin{empheq}{align}\label{kern-K}
K(u,r)  =\frac{a^{d-2\lambda+1}}{u^2-r^2} \, \bigg\{ & \frac{1}{r^2}\,  {}_2 F_1\left(1,\frac{d+2\lambda-3}{2};\frac{d-2\lambda+3}{2};\left(\frac{a}{r}\right)^2\right) 
\\\nonumber
& -  \frac{1}{u^2}\,  {}_2 F_1\left(1,\frac{d+2\lambda-3}{2};\frac{d-2\lambda+3}{2};\left(\frac{a}{u}\right)^2\right) \bigg\},
\end{empheq}
and let
\begin{equation*}
F(r) = \frac{\Gamma((d+2\lambda-3)/2)\Gamma(3-2\lambda)}{2\, \Gamma((d-2\lambda+3)/2)}\, \frac{d}{dr} \int_a^r \frac{Q(t)\,t^{d-2}\,dt}{(r^2-t^2)^{1-\lambda}}, \quad a \leq r \leq b.
\end{equation*}
We can now rewrite integral equation $(\ref{int-9-ch-5})$ as follows,
\begin{empheq}{align}\label{int-10}
G(r) & -  \frac{\Gamma((d+2\lambda-3)/2)\Gamma(3-2\lambda)}{2\, \Gamma((d-2\lambda+3)/2)}\, \int_a^b  G(u)\, K(u,r)\,du \\\nonumber 
& = F_Q \frac{\Gamma((d+2\lambda-3)/2)}{2\pi^{(d-1)/2}\, \Gamma(\lambda)}\, \bigg\{  \frac{d+2\lambda-3}{2} \, r^{d+2\lambda-4}\, \Beta \left(1-\left(\frac{a}{r}\right)^2;\lambda, \frac{d-1}{2}\right) + \frac {a^{d-1}}{r}\, (r^2-a^2)^{\lambda-1} \bigg\} \\\nonumber 
& - F(r), \quad a \leq r \leq b.
\end{empheq}
Integral equation $(\ref{int-10})$ is a Fredholm integral equation of the second kind. We remark that its kernel $K(u,r)$ is symmetric, that is $K(u,r)=K(r,u)$, which can be easily seen from expression $(\ref{kern-K})$.

It remains to mention that the constant $F_Q$ is determined using the fact that $\mu_Q$ is a probability measure, that is its mass is one. We therefore find that
\begin{equation*}
\int_a^b f(t)\, t^{d-2}\, dt = \frac{\Gamma((d-1)/2)}{2\pi^{(d-1)/2}}.
\end{equation*}

\qed

\section{Acknowledgements} 
\noindent This work was done in partial fulfillment of Ph.D. degree at Oklahoma State University under the supervision of Prof. Igor E. Pritsker.

\end{document}